\documentclass[a4paper,fleqn]{cas-sc}

\usepackage[authoryear,longnamesfirst]{natbib}

\usepackage{etoolbox}
\usepackage{lastpage}
\makeatletter

\ExplSyntaxOn
\cs_gset:Npn \__first_footerline: { 
  \group_begin: 
  \small 
  \sffamily 
  \hbox:n{}  % 清空左侧内容
  \group_end: 
}
\ExplSyntaxOff

\def\ps@cas{%
  \let\@oddhead\@empty
  \let\@evenhead\@empty
  \def\@oddfoot{%
    \hbox to \textwidth{%
      \hfill\thepage\ of \pageref{LastPage}\hfill
    }%
    \vskip\footskip
    \hrule height \footrulewidth
  }%
  \let\@evenfoot\@oddfoot
}

\def\ps@plain{%
  \let\@oddhead\@empty
  \let\@evenhead\@empty
  \def\@oddfoot{%
    \hbox to \textwidth{%
      \hfill\thepage\ of \pageref{LastPage}\hfill
    }%
    \vskip\footskip
    \hrule height \footrulewidth
  }%
  \let\@evenfoot\@oddfoot
}

\makeatother
\def\tsc#1{\csdef{#1}{\textsc{\lowercase{#1}}\xspace}}
\tsc{WGM}
\tsc{QE}
\tsc{EP}
\tsc{PMS}
\tsc{BEC}
\tsc{DE}
\numberwithin{equation}{section}

\usepackage{lastpage}
\newtheorem{theorem}{Theorem}
\newtheorem{lemma}[theorem]{Lemma}

\newtheorem{definition}[theorem]{Definition}

\begin{document}
\let\WriteBookmarks\relax
\def\floatpagepagefraction{1}
\def\textpagefraction{.001}
\shorttitle{Singular limits for non-isentropic compressible rotating fluids}
\shortauthors{Yajia Yu et~al.}
%\begin{frontmatter}

\title[mode = title]{Singular limits for non-isentropic compressible rotating fluids}                      
%\tnotemark[1,2]

\author[1]{Yajia Yu}

\author[1]{Chenxi Su}

\author[1]{Ming Lu}
\cormark[1]
\cortext[cor1]{Corresponding author}
\ead{luming@hebtu.edu.cn}

\affiliation[1]{organization={School of Mathematical Sciences, Hebei Normal University}, 
                city={Shijiazhuang},
                postcode={050024}, 
                country={P.R. China}}

\begin{abstract}
In this article, we study the singular limit of non-isentropic compressible rotating fluids. We incorporate the capillary effect into both the $\alpha=1$ and $\alpha=0$ cases, and investigate the Navier-Stokes-Korteweg equations involving the terms of low Mach number, low Rossby number and high Reynolds number. When $\alpha=1$, the dispersion estimate of the acoustic wave equation is derived by Rage's theorem. When $\alpha=0$, we obtain the convergence results by error estimate.
Moreover, we obtain that the three dimensions compressible Navier-Stokes-Korteweg equations converge to the two dimensions incompressible Euler equations.
\end{abstract}


\maketitle

\section{Introduction}

We recall that the non-isentropic compressible Navier-Stokes-Korteweg (Navier- Stokes-Korteweg ) equations 
\begin{equation}\label{1.1}
\begin{cases}
\partial_t\rho+\operatorname{div}(\rho \mathbf u)=0,\\
  \partial_t(\rho \mathbf u)+\operatorname{div}(\rho \mathbf u\otimes{\mathbf u})+\mathbf{e}_3\times (\rho\mathbf u)=\operatorname{div}(S+K),\\
  \partial_t(\rho(\frac{1}{2}\vert\mathbf u\vert^2+e))+\operatorname{div}(\rho \mathbf u (\frac{1}{2}\vert\mathbf u\vert^2+e))=\lambda\Delta\Theta+\operatorname{div}((S+K)\cdot\mathbf u),
\end{cases}
\end{equation}where  $\rho$,  $\mathbf{u}=(\mathbf{u}_1,\mathbf{u}_2,\mathbf{u}_3)^t$,  $e$,  $\Theta$ 
 represent the density, velocity field, internal energy, temperature,
$\lambda\textgreater0$ represent heat conduction,  $\mathbf{e}_3=(0,0,1)$ denotes the axis of rotation, the viscous stress tensor $S$ and the Korteweg tensor $K$ satisfy
\begin{equation*}
    \begin{cases}
     S_{i,j}=2\mu D_{i,j}(\mathbf u)+(\nu\operatorname{div}\mathbf{u}-P)\delta_{i,j},\\
 K_{i,j}=\frac{\kappa}{2}(\Delta\rho^2+\vert\nabla\rho\vert^2)\delta_{i,j}-\kappa\frac{\partial\rho}{\partial x_i}\frac{\partial\rho}{\partial x_j},
    \end{cases}
\end{equation*}
where $D_{i,j}=\frac{1}{2}(\frac{\partial\mathbf{u}_{i,j}}{\partial x_i}+\frac{\partial\mathbf{u}_{i,j}}{\partial x_j})$ is the stress tensor,  $P$ denotes the pressure,  $\mu$ and $\nu$ are the viscosity coefficients satisfying $2\mu+3\nu\textgreater0$ and $\mu\textgreater0$, $\kappa\textgreater0$ is the capillary coefficient. We set $\lambda, \mu$  and $\nu$ to be constants.

Let $\varepsilon\in(0,1]$ is a small covariate, consider the following scale transformations:
\begin{align*}\label{l}
\rho(x,t)=\rho^\varepsilon(x,\varepsilon t), \mathbf u=&\varepsilon\mathbf{u}^\varepsilon(x,\varepsilon t),\Theta(x,t)=\Theta^\varepsilon(x,\varepsilon t),\\
\mu\to \varepsilon\mu&, \nu\to\varepsilon\nu, \lambda\to\varepsilon\lambda.
\end{align*}
Then we obtain the appropriate dimensionless form of $(1.1)$ as follows:
\begin{equation}\label{1.2}
\begin{cases}{}
\partial_{t}\rho^{\varepsilon}+\operatorname{div}\left( \rho^\varepsilon \mathbf{u}^\varepsilon \right)=0,\\
\partial_{t}\left(\rho^\varepsilon \mathbf{u}^\varepsilon\right)+\operatorname{div}\left(\rho^\varepsilon \mathbf {u}^\varepsilon\otimes{\mathbf{u}^\varepsilon}\right)+\dfrac{1}{\text{Ro}}\rho^\varepsilon\left(\mathbf{e}_3\times \mathbf {u}^\varepsilon\right)+\dfrac{1}{\text{Ma}^2}\nabla P^\varepsilon\\
=\dfrac{1}{\text{Re}}\left(\mu+\nu\right)\nabla\operatorname{div}\mathbf{u}^\varepsilon+\dfrac{1}{Re}\mu\Delta\mathbf{u}^\varepsilon+\dfrac{1}{\varepsilon^2}\kappa\rho^\varepsilon
  \nabla\Delta\rho^\varepsilon,\\
 \rho^\varepsilon(e^\varepsilon_t+\mathbf{u}^\varepsilon\cdot\nabla e^\varepsilon)+P^\varepsilon \operatorname{div}\mathbf{u}^\varepsilon-\lambda\Delta\Theta^\varepsilon=\varepsilon^2\Psi(\mathbf{u}^\varepsilon)+\Phi(\rho^\varepsilon,\mathbf{u}^\varepsilon),
\end{cases}
\end{equation}
where  $\text{Ro}$  is the Rossby number,  $\text{Ma}$  is the Mach number and  $\text{Re}$  is the Reynolds number. In this paper, let  $\text{Ma}$  and  $\text{Ro}$  be proportional to  $\varepsilon$  and  $\text{Re}$  be inversely proportional to  $\varepsilon$.  $\Psi(\mathbf{u}^\varepsilon)$  and  $\Phi(\rho^\varepsilon,\mathbf{u}^\varepsilon)$ are denoted as
\begin{align*}\label{l}
\Psi(\mathbf u^\varepsilon)&=\nu(\operatorname{div}\mathbf u^\varepsilon)^2+2\mu D(\mathbf u^\varepsilon):\nabla\mathbf u^\varepsilon,\\
\Phi(\rho^\varepsilon,\mathbf u^\varepsilon)=\kappa(\rho^\varepsilon\Delta\rho^\varepsilon&+\frac{\vert\nabla\rho^\varepsilon\vert^2}{2})\operatorname{div}\mathbf u^\varepsilon-\kappa(\nabla\rho^\varepsilon\otimes\nabla\rho^\varepsilon):\nabla\mathbf u^\varepsilon.
\end{align*}

In this paper, we focus on the compressible fluid obeying the perfect gas relations
\begin{align*}
    P=R \rho \Theta,~e=c_v  \Theta,
\end{align*}
where  $R\textgreater0$,  $c_v\textgreater0$  are the gas constant and the heat capacity at constant volume, respectively. 

For the simplicity,  we will set the specific heat ratio  $\gamma=1+\frac{R}{c_v}$. We denote the density and temperature variations by  $q^\varepsilon$  and  $\theta^\varepsilon$:
\begin{align*}
    \rho^\varepsilon=1+\varepsilon q^\varepsilon, \Theta^\varepsilon=1+\varepsilon \theta^\varepsilon.
\end{align*}
Then equation $(\ref{1.2})$ can be rewritten as
\begin{equation}\label{1.3}
    \begin{cases}{}
\partial_t q^\varepsilon+\dfrac{1}{\varepsilon}\operatorname{div}\mathbf u^\varepsilon+\operatorname{div}\left(q^\varepsilon \mathbf u^\varepsilon\right)=0,\\
  \rho^\varepsilon\left(\mathbf u_t^\varepsilon+\mathbf u^\varepsilon\cdot\nabla \mathbf u^\varepsilon \right)+\dfrac{1}{\varepsilon}\left(\nabla q^\varepsilon+\nabla\theta^\varepsilon\right)+\nabla\left(q^\varepsilon\theta^\varepsilon\right)+\dfrac{1}{\varepsilon}\left(\mathbf{e}_3\times \mathbf u^\varepsilon\right)+\mathbf{e}_3\times \left(q^\varepsilon \mathbf u^\varepsilon\right)\\
 =\varepsilon\mu\Delta\mathbf{u}^\varepsilon+\varepsilon\left(\mu+\nu\right)\nabla\operatorname{div}\mathbf{u}^\varepsilon
  +\dfrac{1}{\varepsilon^{2(1-\alpha)}}\rho^\varepsilon\nabla\Delta\rho^\varepsilon,\\
  \rho^\varepsilon\left(\theta^\varepsilon_t+ \mathbf u^\varepsilon\cdot\nabla\theta^\varepsilon\right)+\left(\rho^\varepsilon\theta^\varepsilon+q^\varepsilon\right)\operatorname{div}\mathbf u^\varepsilon+\dfrac{1}{\varepsilon}\operatorname{div}\mathbf u^\varepsilon-\lambda\Delta\theta^\varepsilon=\varepsilon\Psi\left(\mathbf{u}^\varepsilon\right)+\Phi'\left(\rho^\varepsilon,\mathbf{u}^\varepsilon\right),
 \end{cases}
\end{equation}
where $\Phi'(\rho^\varepsilon,\mathbf{u}^\varepsilon)=\kappa(\rho^\varepsilon\Delta q^\varepsilon+\varepsilon\frac{\vert\nabla q^\varepsilon\vert^2}{2})
  \operatorname{div}\mathbf{u}^\varepsilon-\kappa\varepsilon(\nabla q^\varepsilon\otimes\nabla q^\varepsilon):\nabla\mathbf{u}^\varepsilon$. 

We consider the problem (\ref{1.1})-(\ref{1.3}) in an infinitely extended slab $\Omega$ confined between two parallel planes,
\begin{align*}
    \Omega=R^2\times (0,1).
\end{align*}

The boundary conditions employed herein specify that the velocity field $\mathbf{u}$ adheres to the full-slip constraint.
\begin{align*}
\mathbf u^\varepsilon\cdot{\mathbf n}|_{\partial\Omega}=0,\quad
{\mathbf n}\times \operatorname{curl}\mathbf u^\varepsilon|_{\partial\Omega}=0,\quad
\dfrac{\partial\theta^\varepsilon}{\partial \mathbf n}|_{\partial\Omega}=0,\quad
\nabla\rho^\varepsilon\cdot{\mathbf n}|_{\partial\Omega}=0,
\end{align*}
where $\mathbf{n}=[0,0,\pm1]$ is the outer normal vector.

When the capillary coefficient $\kappa=0$, the Navier-Stokes-Korteweg  equations are reduced to the compressible the Navier-Stokes equations.
Lions and Masmoudi \cite{lions1998incompressible} studied various asymptotic results for the overall weak solution of the compressible Navier-Stokes equations in 1998, and proved that in the low Mach number limit, the fluid flow becomes incompressible. Feireisl and Neasová et al. \cite{feireisl2014inviscid} studied the incompressible limit of the Navier-Stokes equations over a variable region, and proved that the limiting equations of the compressible Navier-Stokes equations is the incompressible Euler equations when the Mach number tends to 0 and the radius of the region is larger than the acoustic velocity. Feireisl and  Gallagher et al. \cite{feireisl2012singular} studied the singular limit problem for the compressible Navier-Stokes equations with rotational effects, and proved that when the \text{Rossby} number and the Mach number both tend to 0, the limiting system of the Navier-Stokes equations with rotational effects is the two-dimensional Navier-Stokes equations.  Feireisl and Novotn\'{y} in  \cite{feireisl2014scale} studied the triple singular limit problem for compressible rotating fluids and proved that the limit equations are modeled as incompressible viscous horizontal motion when the Mach number and the Rossby number both tend to 0 and the Reynolds number tends to infinity, and that the limit equations have a similar structure to those of the two-dimensional Euler equations. Ju and Ou \cite{ju2022low} studied the low Mach number limit of the compressible Navier-Stokes equations with large temperature variations in a three-dimensional bounded region.

When the capillary phenomenon is considered in the viscous compressible Navier-Stokes equations, that is, when the capillary coefficient $\kappa\neq0$, the compressible Navier-Stokes-Korteweg equations are obtained. Van der Waals \cite{van1894thermodynamische} and Korteweg \cite{korteweg1901forme} were the first to study the theory of capillary with diffusive interface. The capillary effect is a phenomenon that occurs when a liquid enters a tiny pore, capillary tube, or tiny channel. When the liquid enters these tiny channels, the liquid will rise or fall due to surface tension. Fanelli \cite{fanelli2016highly} studied the singular limit problem for the isentropic Navier-Stokes-Korteweg equations, considering the incompressible limit and the high rotation limit of the Navier-Stokes-Korteweg equations when the capillary effect is constant or when it vanishes. Kwon and Li \cite{kwon2018incompressible} studied the incompressible limit of the compressible Navier-Stokes-Korteweg  equations, specifically the incompressible limit of the Navier-Stokes-Korteweg equations with a general initial condition on the periodic region $T^3$ and the whole space $R^3$. Hattori and Li \cite{hattori1996existence}considered the solution of local existence and global existence for the three-dimensional compressible Navier-Stokes-Korteweg equations. Hou et al. \cite{hou2018global} considered the three-dimensional non-isentropic compressible Navier-Stokes-Korteweg equations admit a unique global classical solution which has the small initial energy. Hou, Yao and Zhu \cite{hou2017vanishing} investigated that when the capillary coefficient $\kappa$ is sufficiently small, the global smooth solutions converge to the smooth solutions of the three-dimensional compressible non-isentropic NS equations. Sha and Li \cite{sha2019low} studied the low Mach number limit of the three-dimensional nonisentropic compressible Navier-Stokes-Korteweg equations.

In this paper, we study the triple singular limit problem for rotating fluids with capillary effects at low Mach, low Rossby, and high Reynolds numbers. There are some examples of the three singular limits for low Mach, low Rossby and high Reynolds numbers.

The Mach number is the ratio of the velocity of a fluid to the speed of sound, and the speed of sound (i.e. the speed of propagation of sound) has different values at different altitudes, temperatures, and atmospheric densities. In the low Mach number limit, the fluid flow becomes incompressible. Reference \cite{alazard2006low,lions1998incompressible,masmoudi2007examples}.

The Rossby number is a similarity criterion in geohydrodynamics that expresses the ratio of inertia to the Coriolis scale. The Rossby number can be used to characterize the degree of influence of Coriolis forces during planetary rotation. The Rossby number is often used in geophysical phenomena related to the oceans and the Earth's atmosphere. A low Rossby number corresponds to a fast rotation, where a rapidly rotating fluid becomes a two-dimensional plane. Reference \cite{babin1999global,babin2001navier,feireisl2012singular}.

The Reynolds number can be expressed as the ratio of the inertial force to the viscous force of a fluid. The Reynolds number can be used to distinguish between laminar and turbulent flow and to determine the amount of resistance to flow. In the high Reynolds number limit, the viscosity of the fluid is negligible. Reference\cite{clopeau1998vanishing,feireisl2014scale,masmoudi2007examples}.

In this paper, we study the limit problem of the Navier-Stokes equations when $\text{Ma}=\text{Re}=\varepsilon\to0$, $\text{Re}=\frac{1}{\varepsilon}\to\infty$. Fluid flow is expected to degenerate into 2D, inviscid and incompressible.

When $\alpha=1$, we impose the following initial data:
\begin{align}\label{1.4}
\left(q^\varepsilon, \mathbf{u}^\varepsilon,\theta^\varepsilon\right)\vert_{t=0}=\left(q^\varepsilon_0, \mathbf{u}^\varepsilon_0,\theta^\varepsilon_0\right),  \qquad \operatorname{in}\quad\Omega=R^2\times\left(0,1\right),\nonumber\\
(q^\varepsilon_0,\mathbf u^\varepsilon_0,\theta^\varepsilon_0)\operatorname{in} H^3(\Omega),\quad
(q^{\varepsilon}_0,\mathbf u^{\varepsilon}_0,\theta^{\varepsilon}_0)\to(q_0,\mathbf u_0,\theta_0)\quad \operatorname{in}~H^3(\Omega).
\end{align}

If we assume that 
\begin{align*}
 q^\varepsilon\to q, ~~\mathbf{u}^\varepsilon\to \mathbf{u},~~\theta^\varepsilon\to\theta,
\end{align*}
in some sense, then we can derive diagnostic equations
\begin{equation*}
 \begin{cases}
 \mathrm{div}\mathbf u=0,\\
 \mathbf{e}_3\times \mathbf u+\nabla_x\sigma=0,
 \end{cases}
\end{equation*}
where $\sigma=q+\theta$. And the following conditions are met
\begin{align*}
q,~\theta~\operatorname{independent~of} x_3,~~~\sigma=\sigma(x_h),\\
\mathbf u=\mathbf u(x_h),~~~\mathrm{div}_x\mathbf u=\mathrm{div}_h\mathbf u_h=0.
\end{align*}

In the sense of weak solution, the limit equations similar to the 2D Euler equations are obtained
\begin{align}\label{1.5}
  \partial_t(\Delta_h\sigma-\sigma)+\nabla_h^\perp \sigma\cdot\nabla_h(\Delta_h\sigma)=0,
\end{align}
where $\sigma=q+\theta$.

Ou and Li established the following uniform estimates theorem for the strong solution with ill-prepared initial data in \cite{ou2023low}.
\begin{theorem}\label{th1.1}\cite{ou2023low}
 (Uniform estimates for ill-prepared initial data) Let $\Omega\subset R^3$ be a domain with smooth boundary $\partial\Omega$ and $\varepsilon\in(0,1]$ is an arbitrary constant. Suppose that the initial data
$(q^\varepsilon_0,\mathbf u^\varepsilon_0,\theta^\varepsilon_0)\in H^3(\Omega)$ satisfies the assumptions in Lemma \ref{lemma2.1} and that $M_0^\varepsilon\leq D_0,$ where $D_0>0$, is a
constant independent of $\varepsilon$. Assume in addition that there exists a sequence $(q^\varepsilon_0,\mathbf u^\varepsilon_0,\theta^\varepsilon_0)\in H^5(\Omega),$ such
that $(q^{\varepsilon,n}_0,\mathbf u^{\varepsilon,n}_0,\theta^{\varepsilon,n}_0)\to(q^\varepsilon_0,\mathbf u^\varepsilon_0,\theta^\varepsilon_0)$ in $H^3(\Omega)$ as $n\to\infty$ and $(q^{\varepsilon,n}_0,\mathbf u^{\varepsilon,n}_0,\theta^{\varepsilon,n}_0)$ satisfies the
compatibility conditions through order two. If $(q^\varepsilon, \mathbf u^\varepsilon,\theta^\varepsilon)$ is the unique strong solution to (\ref{1.3}) for $\varepsilon$, which is obtained in Lemma \ref{lemma2.1}, then there exist positive constants $\varepsilon_0$, $T$ and $C$, independent of $\varepsilon$, such
that the following uniform estimate holds:
\begin{align*}
M^\varepsilon(t)\leq C,~~\varepsilon\in(0,\varepsilon_0],~~~t\in[0,T].
\end{align*}
\end{theorem}

Based on Theorem \ref{th1.1}, there is the main conclusion in this context.
\begin{theorem}\label{th1.2} Let us take $\alpha=1$:

Let $\varepsilon\in(0,1]$ be fixed and $\Omega=R^2\times (0,1)$. Let
$\sigma_0=q_0+\theta_0=q_0(x_h)+\theta_0(x_h)$ be the unique solution to the elliptic problem, satisfying
\begin{align*}
\Delta_h\sigma_0^\varepsilon+\sigma_0^\varepsilon=\int_{0}^{1}\operatorname{curl}_h \mathbf {u}^\varepsilon_0 dx_3+\int_{0}^{1}\sigma_0^\varepsilon dx_3.
\end{align*}
The initial value problem (\ref{1.3}) has a unique strong solution $(q^\varepsilon,\mathbf u^\varepsilon,\theta^\varepsilon)$, which satisfies as $\varepsilon\to0$ that
\begin{align*}
q^\varepsilon\stackrel{w}{\longrightarrow} q~~~&in~L^2(0,T;H^2(\Omega));\qquad q^\varepsilon\stackrel{w*}{\longrightarrow} q~~~in~L^\infty(0,T;H^1(\Omega)),\\
\mathbf u^\varepsilon\stackrel{w}{\longrightarrow} \mathbf u~~~&in~L^2(0,T;H^3(\Omega));\qquad \mathbf u^\varepsilon\stackrel{w*}{\longrightarrow} \mathbf u~~~in~L^\infty(0,T;H^2(\Omega)),\\
\theta^\varepsilon\stackrel{w}{\longrightarrow} \theta~~~&in~L^2(0,T;H^3(\Omega));\qquad \theta^\varepsilon\stackrel{w*}{\longrightarrow} \theta~~~in~L^\infty(0,T;H^1(\Omega)),
\end{align*}
and
\begin{align*}
&\rho^\varepsilon=1+\varepsilon q^\varepsilon\to 1~~~in~L^2(0,T;H^2(\Omega))\cap L^\infty(0,T;H^1(\Omega)),\\
&\Theta^\varepsilon=1+\varepsilon\theta^\varepsilon\to 1~~~in~L^2(0,T;H^2(\Omega))\cap L^\infty(0,T;H^1(\Omega)),\\
&\mathbf u^\varepsilon\to \mathbf u,~~in~L^2((0,T)\times K;R^3), for~any~compact~K\subset\Omega,
\end{align*}
where $(q,\mathbf{u},\theta)$ is the solution of~(\ref{1.4}).
\end{theorem} 

When $\alpha=0$,  we impose the following initial:
\begin{align*}
\left(q^\varepsilon, \mathbf{u}^\varepsilon,\theta^\varepsilon\right)\vert_{t=0}=\left(q^\varepsilon_0, \mathbf{u}^\varepsilon_0,\theta^\varepsilon_0\right),  \qquad \operatorname{in}\quad\Omega=R^2\times\left(0,1\right),\\
\mathbf{w}_0(x)\in H^{s+2}(\Omega),\quad
\Vert q_0^\varepsilon\Vert_{H^{s+1}(\Omega)}+\Vert \mathbf{u}_0^\varepsilon-\mathbf{w}_0\Vert_{H^{s}(\Omega)}+\Vert \theta_0^\varepsilon\Vert_{H^{s}(\Omega)}=O(\varepsilon).
\end{align*}

If we assume that 
\begin{align*}
    q^\varepsilon\to0, ~~\mathbf{u}^\varepsilon\to\omega,~~\theta^\varepsilon\to0,
\end{align*}
in some sense, we can derive diagnostic equations
\begin{equation*}
 \begin{cases}
 \mathrm{div}\mathbf u=0,\\
 \mathbf{e}_3\times \mathbf u+\nabla_x r=0,
 \end{cases}
\end{equation*}
where $r=(\mathrm{I}-\Delta)q+\theta$. And the following conditions are met
\begin{align*}
q,~\theta~\operatorname{independent~of} x_3,~~~r=r(x_h),\\
\mathbf u=\mathbf u(x_h),~~~\mathrm{div}_x\mathbf u=\mathrm{div}_h\mathbf u_h=0.
\end{align*}
And the limit equations are
\begin{equation}\label{1.6}
  \begin{cases}
  \operatorname{div}\mathbf{w}=\operatorname{div}\mathbf{w}_h=0,\\
  \mathbf{w}_t+\mathbf{w}\cdot\nabla_h\mathbf{w}+\nabla_h\pi=0,
  \end{cases}
\end{equation}
where the limits of $\mathbf{u}^\varepsilon$ and $\frac{1}{\varepsilon}[(1+\varepsilon q^\varepsilon)\nabla\theta^\varepsilon+(1+\varepsilon\theta^\varepsilon)\nabla q^\varepsilon-(1+\varepsilon q^\varepsilon)\nabla\Delta q^\varepsilon-(1+\varepsilon q^\varepsilon)\nabla((I-\Delta)q^\varepsilon+\theta^\varepsilon)]$ by $\mathbf w$ and $\nabla\pi$.

The following theorem is another major result in this paper.

\begin{theorem}\label{th1.3} Let us take $\alpha=0$:
 
Let $s\textgreater\frac{3}{2}+2$, then there exists a constant $\varepsilon_0\textgreater0$ such that for all $\varepsilon\in(0,\varepsilon_0]$. The Navier-Stokes-Korteweg equations (\ref{1.3}) have a unique solution $(q^\varepsilon,\mathbf{u}^\varepsilon,\theta^\varepsilon)(t,x)$ that satisfies $q^\varepsilon\in C([0,T^\ast],H^{s+1}(\Omega))$ and $(\mathbf{u}^\varepsilon,\theta^\varepsilon)\in C([0,T^\ast],H^s(\Omega))$. And there exists a positive constant $K$ that is independent of $\varepsilon$ such that for all $\varepsilon\in(0,\varepsilon_0],$
\begin{align}\label{1.7}
\operatorname{sup}\{\Vert q^\varepsilon-\varepsilon^2\pi\Vert_{H^{s+1}(\Omega)}+\Vert\mathbf{u}^\varepsilon-\mathbf{w}
\Vert_{H^s(\Omega)}+\Vert\theta^\varepsilon-\varepsilon\pi\Vert_{H^s(\Omega)}\}\leq K\varepsilon,~~~t\in[0,T^\ast],
\end{align}
where $(\mathbf{w},\pi)\in C([0,T^\ast]$, $H^{s+2}(\Omega))\cap C^1([0,T^\ast],H^s(\Omega))$ is a smooth solution of the incompressible Navier-Stokes equations (\ref{1.6}).
\end{theorem}

The paper is organized as follows: In Section \ref{s2}, we present the fundamental theorems employed in our analysis. In Section \ref{s3}, we first establish the priori estimates without capillary effects (i.e. when $\alpha=1$), then derive the corresponding dispersion estimates, and finally investigate the limiting systems. In Section \ref{s4}, we analyze the error estimate incorporating capillary phenomena (i.e. when $\alpha=0$).
\section{Preliminaries}\label{s2}
This section presents the theorems that are required for the subsequent steps. First, we present the local-in-time existence theory for (\ref{1.3}) with the initial data $(q^\varepsilon_0,\mathbf u^\varepsilon_0,\theta^\varepsilon_0)$.
\begin{lemma}\cite{korteweg1901forme}\label{lemma2.1}
     Let $\varepsilon\in(0,1]$ be fixed and $\Omega\in R^3$ be a singly connected region with smooth boundary $\partial\Omega$. Initial value condition satisfied $(q^\varepsilon_0,\mathbf u^\varepsilon_0,\theta^\varepsilon_0)\in H^{s+1}(\Omega)\times H^s(\Omega)\times H^s(\Omega)$, assume that the following compatibility conditions are satisfied
\begin{align*}
\mathbf u^\varepsilon_0\cdot{\mathbf n}=\mathbf u^\varepsilon_t(0)\cdot{\mathbf n}=0,~~{\mathbf n}\times \operatorname{curl}\mathbf u^\varepsilon_0={\mathbf n}\times \operatorname{curl}\mathbf u^\varepsilon_t(0)=0,~~in~\partial\Omega,
\end{align*}
\begin{align*}
\frac{\partial\theta_0^\varepsilon}{\partial{\mathbf n}}=\frac{\partial\theta_t^\varepsilon(0)}{\partial{\mathbf n}}=0,~~in~\partial\Omega.
\end{align*}
Then there exists a positive constant $T$ such that the system of equations (\ref{1.3}) with initial value condition $(q_0^\varepsilon, \mathbf{u}_0^\varepsilon, \theta_0^\varepsilon)(x)$ has a unique solution $(q^\varepsilon,\mathbf{u}^\varepsilon,\theta^\varepsilon)$, satisfies
\begin{align*}
q^\varepsilon\in C([0,T],H^{s+1}(\Omega))~and~(\mathbf{u}^\varepsilon,\theta^\varepsilon)\in C([0,T],H^s(\Omega)).
\end{align*}
\end{lemma}

Next, we give the lemma for the 2D Navier-Stokes equations on the existence of a unique regular solution in local time.
\begin{lemma}\cite{kato1984nonlinear}\label{lemma2.2}
Assume that $s\geq3$. Let $\operatorname{div}\mathbf{w}_0\in H^s(R^2;R^2)$ and satisfy $\mathbf{\omega}_0=0$. It is well known that the Euler system (\ref{1.6}) supplemented with the initial data $\mathbf{w}(0)=\mathbf{w}_0$, admits a regular solution $(\mathbf{w},\pi)$, unique in the class
\begin{align*}
\mathbf{w}\in C([0,T];H^s(R^2;R^2)),~~\partial_t\mathbf{w}\in C([0,T];H^{s-1}(R^2;R^2)),~~\pi\in C([0,T];H^s(R^2)).
\end{align*}
\end{lemma}

Sha and Li in \cite{sha2019low} define the maximal existence interval of the solution, and show the $\lim\limits_{\varepsilon \to 0}T_\varepsilon >0$.
\begin{lemma}\cite{sha2019low}\label{lemma2.3}
Define
\begin{align*}
T_\varepsilon:=\operatorname{sup}\{T\textgreater0:q^\varepsilon(t,x)\in C([0,T],H^{s+1}(\Omega)),(\mathbf{u}^\varepsilon,\theta^\varepsilon)(t,x)\in C([0,T],H^s(\Omega));\\
c_1\leq1+\varepsilon q^\varepsilon(t,x)\leq c_2,~~c_1\leq1+\varepsilon\theta^\varepsilon(t,x)\leq c_2~~~\forall(x,t)\in\Omega\times[0,T] \}.
\end{align*}
Namely, $[0,T_\varepsilon)$ is the maximal time interval of $H^{s+1}(\Omega)\times H^{s}(\Omega)\times H^{s}(\Omega)$ existence.

Under the conditions of Theorem \ref{th1.3}, there exists a constant $\varepsilon_0=\varepsilon_0(T^\ast)$ such that for all
$\varepsilon\leq\varepsilon_0$, it holds that
\begin{align*}
T_\varepsilon\textgreater T^\ast.
\end{align*}
\end{lemma}

In the study of acoustic analysis, we need to invoke the RAGE theorem.
\begin{lemma}[RAGE theorem]\cite{feireisl2012singular}\label{lemma2.4}
     Let  $H$  be a Hilbert space,  $A: \mathcal{D}(A) \subset H \rightarrow H$  a selfadjoint operator, $C: H \rightarrow H$  a compact operator, and  $P_{c}$  the orthogonal projection onto $H_{c}$, where
 \begin{align}
H & = H_{c} \oplus \operatorname{cl}_{H}\{\operatorname{span}\{w \in H \mid w \text { an eigenvector of } A\}\}.\nonumber 
\end{align}
Then
\begin{align}
\left\|\frac{1}{\tau} \int_{0}^{\tau} \exp (-\operatorname{i} t A) C P_{c} \exp (\operatorname{i} t A) \operatorname{d} t\right\|_{\mathcal{L}(H)} \rightarrow 0 \text { for } \tau \rightarrow \infty.\nonumber 
\end{align}
\end{lemma}

Finally, we list some inequalities, which play an important role in uniform bounds and error estimate.
\begin{lemma}\cite{sha2019low}\label{lemma2.5}
Let $s\textgreater\frac{3}{2}+1$ be an integer and $\alpha=\alpha(\alpha_1,\alpha_2,\alpha_3)$ be a multi-index such that $|\alpha|\leq s$.
Then there exists a positive constant $C$ such that for any $f\in H^s (\Omega)$, $g\in H^{s-1} (\Omega)$, $[D^\alpha,f]g\in L^2(\Omega)$
and it holds that
\begin{align*}
\vert[D^\alpha,f]g\Vert\leq C(\Vert\nabla f\Vert_{L^\infty}\Vert g\Vert_{H^{s-1}}+\Vert f\Vert_{H^s}\Vert g\Vert_{L^\infty}).
\end{align*}
\end{lemma}
\begin{lemma}\cite{BourguignonBrezis1974}\label{lemma2.6}
Let $\Omega$ be a bounded domain in $R^N$ with smooth boundary $\partial\Omega$
and its unit outward normal $\mathbf{n}$. Then there exists a constant $C\textgreater 0$ independent of $\mathbf{u}$, such that
\begin{equation*}
\Vert \mathbf u\Vert_{H^s(\Omega)}\leq C(\Vert \mathrm{div}\mathbf u\Vert_{H^{s-1}(\Omega)}+\Vert \mathrm{curl}\mathbf u\Vert_{H^{s-1}(\Omega)}+\Vert \mathbf u\cdot\mathbf n\Vert_{H^{s-\frac{1}{2}}(\partial\Omega)}+\Vert \mathbf u\Vert_{H^{s-1}(\Omega)}),
\end{equation*}
for any $\mathbf{u}\in[H^s(\Omega)]^N, s\ge1$.
\end{lemma}
\begin{lemma}\cite{ou2023low}\label{lemma2.7}
Let $\Omega$ be a bounded domain in $R^3$ with smooth boundary $\partial\Omega$ and its
unit outward normal $\mathbf{n}$. Suppose that $\mathbf{n}\times{\mathrm{curl}\mathbf{u}}$=0 on $\partial\Omega$ for any vector $\mathbf{u}\in H^{k+2}(\Omega)$, $k=0,1,2,\cdots$.
Then there exists a constant $C=C(k,\Omega)>0$ independent of $\mathbf{u}$, such that
\begin{equation*}
\Vert \mathrm{curl}\mathrm{curl} \boldsymbol u\cdot\boldsymbol n\Vert_{H^{k+\frac{1}{2}}(\partial\Omega)}\leq C\Vert \boldsymbol u\Vert_{H^{k+2}(\Omega)}.
\end{equation*}
\end{lemma}

\section{Limit problem of capillary effect when $\alpha=1$}\label{s3}
This section proves the Theorem \ref{th1.2} of this thesis.
\subsection{Uniform bounds}\label{s3.1}

\begin{definition}\label{def1.1}
We define the energy inequality:
\begin{align*}
M^\varepsilon(t):=\mathop{\operatorname{esssup}}\limits_{0\leq s\leq t}\{&\Vert(\mathbf u^\varepsilon,\rho^\varepsilon,\varepsilon q^\varepsilon,\varepsilon\theta^\varepsilon)\Vert^2_{H^2(\Omega)}(s)+\Vert(q^\varepsilon,\theta^\varepsilon)\Vert^2_{H^1(\Omega)}(s)+\Vert\varepsilon(q_t^\varepsilon,\mathbf u_t^\varepsilon,\theta_t^\varepsilon)
\Vert^2_{H^1(\Omega)}(s)\\
+&\Vert\varepsilon^2(q_{tt}^\varepsilon,\mathbf u_{tt}^\varepsilon,\theta_{tt}^\varepsilon)\Vert^2_{L^2(\Omega)}(s)+\Vert(\varepsilon \mathbf u^\varepsilon,\varepsilon^2 q^\varepsilon,\varepsilon^2\theta^\varepsilon)\Vert^2_{H^3(\Omega)}(s)
+\Vert(\rho^\varepsilon,(\rho^\varepsilon)^{-1})\Vert^2_{L^\infty(\Omega)}(s)\}\\
+\int_{0}^{t}\{&\Vert(q^\varepsilon,\theta^\varepsilon)\Vert^2_{H^2(\Omega)}+\Vert(\varepsilon q^\varepsilon,\varepsilon\mathbf u^\varepsilon,\varepsilon\theta^\varepsilon)\Vert^2_{H^
3(\Omega)}+\Vert\varepsilon(q_t^\varepsilon,\mathbf u_t^\varepsilon,\theta_t^\varepsilon)\Vert^2_{H^2(\Omega)}\\
+&\varepsilon\Vert\varepsilon^2\mathbf u_{tt}^\varepsilon\Vert^2_{H^1(\Omega)}+\Vert\varepsilon^2\theta_{tt}^\varepsilon\Vert^2_{H^1(\Omega)}+\varepsilon\Vert\varepsilon^2\mathbf u^\varepsilon\Vert^2_{H^4(\Omega)}+\Vert(\varepsilon^2q^\varepsilon,\varepsilon^2\theta^\varepsilon)\Vert^2_{H^4(\Omega)}\}ds.
\end{align*}
The initial weighted energy is
\begin{align*}
M_0^\varepsilon:=M^\varepsilon(t=0).
\end{align*}
\end{definition}

Using $(\ref{1.3})_1$ multiplied by $\rho^{-k}$, referring to \cite{ou2022incompressible} and combining with the Gronwall inequality, we get
\begin{align}\label{3.1}
\Vert(\rho,\rho^{-1})\Vert^2_{L^\infty}(t)+\Vert\rho\Vert^2_{H^2}(t)\leq C_0(M_0)\operatorname{exp}{t^{\frac{1}{2}}C_1(M(t))}.    
\end{align}
Integrate by multiplying  $(\ref{1.3})_1$,  $(\ref{1.3})_2$ and  $(\ref{1.3})_3$ respectively by $q$, $\mathbf{u}$,  $\theta$, and using Lemma \ref{lemma2.6}, we get
\begin{align}\label{3.2}
\Vert(q,\mathbf u,\theta)\Vert^2_{L^2}(t)+\varepsilon\Vert \mathbf u\Vert^2_{L^2_t{(H^1)}}+\Vert\nabla\theta\Vert^2_{L^2_t{(L^2)}}\leq C_0(M_0)\operatorname{exp}\{(t^\frac{1}{2}+\varepsilon)C_1(M(t))\}.
\end{align}
Next integrate by multiplying $\nabla(\ref{1.3})_1$, $\nabla(\ref{1.3})_3$ by $\nabla q$, $\nabla\theta$, and take the inner product of $\left<\ref{1.3})_2, -\nabla\operatorname{div}\mathbf u\right>$, according to $\mathbf{u}\cdot\mathbf{n}\vert_{\partial\Omega}=0$, we get
\begin{align}\label{3.3}
\Vert(\nabla q,\operatorname{div}\mathbf u,\nabla\theta)\Vert^2_{L^2}(t)+\varepsilon\Vert\nabla \operatorname{div}\mathbf u\Vert^2_{L^2_t(L^2)}+\Vert\Delta\theta\Vert^2_{L^2_t(L^2)}
\leq C_0(M_0)\operatorname{exp}\{(t^\frac{1}{2}+\varepsilon)C_1(M(t))\}.
\end{align}
In the calculation of (\ref{3.3}), when term $\int_{0}^{t}\int_{\Omega}(\boldsymbol{e}_3\times \boldsymbol u)\cdot \nabla\mathrm{div}\boldsymbol udxds=\int_{0}^{t}\int_{\Omega}\nabla^{\bot}_h\boldsymbol u_h\cdot \mathrm{div}\boldsymbol u dxds$ appears, we handle it as follows: multiply $(\ref{1.3})_3$ by $-\nabla^{\bot}_h\boldsymbol u_h$ to take the integral over $\Omega\times(0,t)$, we get
\begin{align*}
-\frac{1}{\varepsilon}\int_{0}^{t}\int_{\Omega}\nabla^{\bot}_h\boldsymbol u_h\cdot \mathrm{div}\boldsymbol u dxds\leq C_0(M_0)\mathrm{exp}\{(t^{\frac{1}{2}}+\varepsilon)C_1(M(t))\}.
\end{align*}
By combining $(\ref{3.1})$, $(\ref{3.2})$, $(\ref{3.3})$, we obtain estimate for $\rho$  and  $\left(q,\mathbf{u},\theta\right)$.

Differentiating  $(\ref{1.3})$  with respect to t, we deduce that
\begin{equation}\label{3.4}
\begin{cases}{}
 q_{tt}+\dfrac{1}{\varepsilon}\operatorname{div}\mathbf u_t=-(\mathbf u_t\cdot\nabla q+\mathbf u\cdot\nabla q_t+q_t\cdot \operatorname{div}\mathbf u+q\operatorname{div}\mathbf u_t)~~~~~~~~~~~~~~~~~~~~~~~~~~~~~~~~~~~~~(t,x)\in(0,T]\times\Omega,\\
  \rho(\mathbf u_{tt}+\mathbf u\cdot\nabla \mathbf u_t)+\dfrac{1}{\varepsilon}\nabla q_t+\dfrac{1}{\varepsilon}\nabla \theta_t-\varepsilon\mu\Delta \mathbf u_t-\varepsilon(\mu+\nu)\nabla \operatorname{div}\mathbf u_t+\dfrac{1}{\varepsilon}(\mathbf{e}_3\times \mathbf u_t)\\
  =-\rho_t(\mathbf u_t+\mathbf u\cdot\nabla \mathbf u)
  -\rho \mathbf u_t\cdot\nabla \mathbf u-\nabla(q\theta)_t-\mathbf{e}_3\times(q\mathbf u)_t+(\rho\nabla\Delta\rho)_t~~~~~~~~~~~~~~~~~~~~~~~~~~~(t,x)\in(0,T]\times\Omega,\\
  (\rho \mathbf u)_t\cdot\nabla\theta+\dfrac{1}{\varepsilon}\operatorname{div}\mathbf u_t+\rho \mathbf u\cdot\nabla\theta_t+\rho_t\theta_t+\rho\theta_{tt}+\operatorname{div}\mathbf u_t(\rho\theta+q)+\operatorname{div}\mathbf u(\rho_t\theta+q_t+\rho\theta_t)\\
  =\lambda\Delta\theta_t+[\varepsilon\Psi(\mathbf{u}^\varepsilon)+\Phi'(\rho^\varepsilon,\mathbf{u}^\varepsilon)]_t~~~~~~~~~~~~~~~~~~~~~~~~~~~~~~~~~~~~~~~~~~~~~~~~~~~~~~~~~~~~~~~~~~~~~~~~~~~~~(t,x)\in(0,T]\times\Omega,\\
  \mathbf u_t\cdot{\mathbf n}=0,~~ {\mathbf n}\times \operatorname{curl}\mathbf u_t=0~~~~~~~~~~~~~~~~~~~~~~~~~~~~~~~~~~~~~~~~~~~~~~~~~~~~~~~~~~~~~~~~~~~~~~~~~~~~~~~~~~~~~~~(0,T]\times\partial\Omega,\\
  (q_t,\mathbf u_t)(0,x)=(q_t(0),\mathbf u_t(0))~~~~~~~~~~~~~~~~~~~~~~~~~~~~~~~~~~~~~~~~~~~~~~~~~~~~~~~~~~~~~~~~~~~~~~~~~~~~~~~~~~~~~~x\in\Omega.
\end{cases}
\end{equation}
where 
\begin{align*}
q_t(0):&=-\dfrac{1}{\varepsilon}\operatorname{div}\mathbf u_0-q_0\operatorname{div}\mathbf u_0-\mathbf u_0\cdot\nabla q_0,\\
\mathbf u_t(0):&=[\varepsilon\mu\Delta \mathbf u_0+\varepsilon(\mu+\nu)\nabla \operatorname{div}\mathbf u_0+\rho_0\nabla\Delta\rho_0-\nabla(q_0\theta_0)-\dfrac{1}{\varepsilon}\nabla q_0-\dfrac{1}{\varepsilon}\nabla\theta_0-\dfrac{1}{\varepsilon}(\mathbf{e}_3\times\mathbf u_0)\\
&-\mathbf{e}_3\times(q_0\mathbf u_0)]
/\rho_0-\mathbf u_0\cdot\nabla \mathbf u_0,\\
\theta_t(0):&=[\lambda\Delta\theta_0+\varepsilon\Psi(\mathbf{u}_0^\varepsilon)+\Phi'(\rho_0^\varepsilon,\mathbf{u}_0^\varepsilon)
-\dfrac{1}{\varepsilon}\operatorname{div}\mathbf u_0-(\rho_0\theta_0+q_0)\operatorname{div}\mathbf u_0]/\rho_0-\mathbf u_0\cdot\nabla\theta_0.
\end{align*}
Integrate by multiplying  $(\ref{3.4})_1$,  $(\ref{3.4})_2$ and  $(\ref{3.4})_3$ respectively by $\varepsilon^2 q_t$, $\varepsilon^2\mathbf{u}_t$,  $\varepsilon^2\theta_t$ in $\Omega\times (0,t)$ and integrating by parts, and in conjunction with Lemma \ref{lemma2.6}, we get
\begin{align}\label{3.5}
\Vert\varepsilon(q_t,\mathbf u_t,\theta_t)\Vert^2_{L^2}(t)+\varepsilon\Vert\varepsilon \mathbf u_t\Vert^2_{L^2_t(H^1)}+\lambda\Vert\varepsilon\nabla\theta_t\Vert^2_{L^2_t(L^2)}\leq C_0(M_0)\operatorname{exp}\{(t^{\frac{1}{2}}+\varepsilon)C_1(M(t))\},
\end{align}
Next, we take the inner product of  $\nabla(\ref{3.4})_1$ and  $\varepsilon^2\nabla q_t$, $(\ref{3.4})_2$ and $-\varepsilon^2\nabla\operatorname{div}\mathbf{u}_t$,  $\nabla(\ref{3.4})_3$ and $\varepsilon^2\nabla\theta_t$ in  $\Omega\times (0,t)$  to get
\begin{align}\label{3.6}
&\Vert\varepsilon(\nabla q_t,\operatorname{div}\mathbf u_t,\nabla\theta_t)\Vert^2_{L^2}(t)+\varepsilon\Vert\varepsilon\nabla \operatorname{div}\mathbf u_t\Vert^2_{L^2_t(L^2)}+\lambda\Vert\varepsilon\Delta\theta_t\Vert^2_{L^2_t(L^2)}\leq C_0(M_0)\operatorname{exp}\{(t^\frac{1}{2}+\varepsilon)C_1(M(t))\}.
\end{align}
By combining $(\ref{3.5})$, $(\ref{3.6})$, we obtain estimate for $\rho$  and  $\left(\varepsilon q_t,\varepsilon\mathbf{u}_t,\varepsilon\theta_t\right)$.

Multiplying  $(\ref{1.3})_2$ by  $\rho^{-1}(\varepsilon q)$
\begin{align}\label{3.7}
&\varepsilon q(\mathbf u_t+\mathbf u\cdot\nabla \mathbf u)+\frac{1}{\varepsilon}\rho^{-1}(\varepsilon q)(\nabla q+\nabla\theta)+\rho^{-1}(\varepsilon q)\nabla(q\theta)+\frac{1}{\varepsilon}\rho^{-1}(\varepsilon q)(\mathbf{e}_3\times \mathbf u)\nonumber\\
&=\rho^{-1}(\varepsilon q)(\varepsilon(2\mu+\nu)\nabla \operatorname{div}\mathbf u-\varepsilon\mu \operatorname{curl}\operatorname{curl}\mathbf u)-\rho^{-1}(\varepsilon q)\mathbf{e}_3\times(q\mathbf u)+\rho^{-1}(\varepsilon q)(\rho\nabla\Delta\rho).
\end{align}
Combining  $\operatorname{curl}(\ref{3.7})$  and  $(\ref{1.3})$, we obtain
\begin{equation}\label{3.8}
\begin{cases}{}
\varepsilon q\omega_t-\mu\rho^{-1}(\varepsilon^2 q)\Delta\omega=f_1=A_1+A_2+A_3+A_4+A_5~~~~~~~~~~~~~~~~~~~~~~~~~~~~~~~~(t,x)\in(0,T]\times\Omega,\\
  \omega\times{\mathbf n}=0~~~~~~~~~~~~~~~~~~~~~~~~~~~~~~~~~~~~~~~~~~~~~~~~~~~~~~~~~~~~~~~~~~~~~~~~~~~~~~~~~~~~~~~~~~~~~~~~~~~~~~~~~~(0,T]\times\partial\Omega,\\
  \omega(0,x)=\omega_0:=\operatorname{curl}\mathbf u_0~~~~~~~~~~~~~~~~~~~~~~~~~~~~~~~~~~~~~~~~~~~~~~~~~~~~~~~~~~~~~~~~~~~~~~~~~~~~~~~~~~~~~~x\in\Omega.
\end{cases}
\end{equation}
where
\begin{equation*}
\begin{cases}
A_1=-\nabla(\varepsilon q)\times \mathbf u_t-\nabla(\varepsilon q)\times(\mathbf u\cdot\nabla \mathbf u)-\varepsilon q\nabla\times(\mathbf u\cdot\nabla \mathbf u),\\
 A_2=-\nabla(\rho^{-1}q)\times(\nabla q+\nabla\theta)-\nabla(\rho^{-1}\varepsilon q)\times\nabla(q\theta),\\
 A_3=-\nabla(\rho^{-1}q)\times(\mathbf{e}_3\times \mathbf u)-\rho^{-1}q\nabla\times(\mathbf{e}_3\times \mathbf u)-\nabla(\rho^{-1}q\varepsilon)\times(\mathbf{e}_3\times q\mathbf u)-(\rho^{-1}q\varepsilon)\nabla\times(\mathbf{e}_3\times q\mathbf u),\\
 A_4=\nabla(\varepsilon^2\rho^{-1} q)\times(-\mu \operatorname{curl}\operatorname{curl}\mathbf u+(2\mu+\nu)\nabla \operatorname{div}\mathbf u),\\
 A_5=\nabla(\varepsilon q)\times(\nabla\Delta\rho).
\end{cases}
\end{equation*}
 And the facts that $\omega=\operatorname{curl}\mathbf u=(\partial_{x_2}\mathbf u_3-\partial_{x_3}\mathbf u_2,\partial_{x_3}\mathbf u_1-\partial_{x_1}\mathbf u_3,\partial_{x_1}\mathbf u_2-\partial_{x_2}\mathbf u_1)^t$, $\Delta \operatorname{curl}=-\operatorname{curl}\operatorname{curl}\operatorname{curl}$, $\operatorname{curl}\nabla=0$ are used in deriving $(\ref{3.8})$.\\
Additionally, $\operatorname{curl}(\ref{1.3})_2$ is multiplied by  $\varepsilon q$  to obtain
\begin{align}\label{3.9}
&\varepsilon q\omega_t+\varepsilon q\nabla\times(\mathbf u\cdot\nabla\mathbf u)+\varepsilon q\nabla \times[\varepsilon q\times(\mathbf {u}_t+\mathbf u\cdot\nabla \mathbf u)]+q\nabla\times(\mathbf{e}_3\times \mathbf u)\nonumber+\varepsilon q\nabla\times(\mathbf{e}_3\times(q\mathbf u))\\
&=\varepsilon^2\mu q\Delta\operatorname{curl} \mathbf{u}+\varepsilon\nabla\times(\rho\nabla\Delta\rho).
\end{align}
By combining Lemma \ref{lemma2.6}, Lemma \ref{lemma2.7}, $(\ref{3.8})$, $(\ref{3.9})$ and boundary condition  $\operatorname{curl}\mathbf{u}\times\mathbf{n}\vert_{\partial\Omega}=0$, we obtain
\begin{align}\label{3.10}
\Vert \operatorname{curl}\mathbf u\Vert^2_{H^1}(t)+\varepsilon\Vert \operatorname{curl}\mathbf u\Vert^2_{L^2_t(H^2)}\leq C_0(M_0)\operatorname{exp}\{t^\frac{1}{2}C_1(M(t))\}.
\end{align}
Next, we apply to  $\operatorname{curl}\partial_t(\ref{3.8})$  to derive
\begin{equation}\label{3.11}
\begin{cases}
\varepsilon q\omega_{tt}-\mu\operatorname{curl}(\rho^{-1}(\varepsilon^2 q)\operatorname{curl}\omega_t)=f_2=B_1+B_2+B_3+B_4~~~~~~~~~~~~~~~~~~~~~~~~(t,x)\in(0,T]\times\Omega,\\
  \omega\times{\mathbf n}=0~~~~~~~~~~~~~~~~~~~~~~~~~~~~~~~~~~~~~~~~~~~~~~~~~~~~~~~~~~~~~~~~~~~~~~~~~~~~~~~~~~~~~~~~~~~~~~~~~~~~~~~~(0,T]\times\partial\Omega,\\
  \omega_t(0,x):=\operatorname{curl}\mathbf u_t(0)~~~~~~~~~~~~~~~~~~~~~~~~~~~~~~~~~~~~~~~~~~~~~~~~~~~~~~~~~~~~~~~~~~~~~~~~~~~~~~~~~~~~~~~~x\in\Omega.
\end{cases}
\end{equation}
where
\begin{equation*}
\begin{cases}
\omega_t(0)=\operatorname{curl}[\rho_0^{-1}(\varepsilon(2\mu+\nu)\Delta \mathbf u_0+\varepsilon\mu\nabla \operatorname{div}\mathbf u_0-\frac{1}{\varepsilon}\mathbf{e}_3\times \mathbf u-\mathbf{e}_3\times(q\mathbf u))]-\operatorname{curl}(\mathbf u_0\cdot\nabla \mathbf u_0),\\
B_1=-\nabla\times(\varepsilon q_t \mathbf u_t)-\nabla\times[\varepsilon q_t\cdot (\mathbf u\cdot\nabla \mathbf u)]-\nabla(\varepsilon q)\times\mathbf u_{tt}-\nabla\times[\varepsilon q(\mathbf u\cdot\nabla \mathbf u)_t],\\
B_2=-\nabla(\rho^{-1}q)_t\times(\nabla q+\nabla\theta)-\nabla(\rho^{-1}q)\times\nabla( q_t+\theta_t)-\varepsilon\nabla\times(\rho^{-1} q)_t\cdot\nabla(q\theta)-\varepsilon\nabla\times(\rho^{-1}q)\cdot\nabla(q\theta)_t,\\
B_3=-\nabla\times[(\rho^{-1}q)_t\cdot(\mathbf{e}_3\times \mathbf u)]-\nabla\times[(\rho^{-1}q)\cdot(\mathbf{e}_3\times \mathbf u)_t]-\varepsilon\nabla\times[(\rho^{-1}q)_t(\mathbf{e}_3\times q\mathbf u)]\\
~~~~~~~~~~-\varepsilon\nabla\times[(\rho^{-1}q)\cdot(\mathbf{e}_3\times q\mathbf u)_t],\\ B_4=\varepsilon^2\nabla(\rho^{-1}q)_t\cdot[(2\mu+\nu)\nabla\operatorname{div}\mathbf u-\mu \operatorname{curl}\operatorname{curl}\mathbf u]+\varepsilon^2\mu(\rho^{-1}q)_t\Delta\operatorname{curl}\mathbf{u}+\varepsilon^2\nabla(\rho^{-1}q)(\mu+\nu)\nabla\operatorname{div}\mathbf{u}_t\\
~~~~~~~~~~+\varepsilon\nabla q_t\times(\nabla\Delta\rho)+\varepsilon\nabla q\times(\nabla\Delta\rho_t).
\end{cases}
\end{equation*}
We multiply  $(\ref{3.11})_1$  by  $\varepsilon^2\omega_t$  and integrate the resulting equality in  $\Omega\times (0,t)$, we obtain
\begin{align}\label{3.12}
\Vert \varepsilon\operatorname{curl}\mathbf u_t\Vert^2_{L^2}(t)+\varepsilon\Vert \varepsilon \operatorname{curl}\mathbf u_t\Vert^2_{L^2_t(H^1)}\leq C_0(M_0)\operatorname{exp}\{(t^\frac{1}{2}+\varepsilon)C_1(M(t))\}.
\end{align}
Then using  $\varepsilon\nabla(\ref{3.9})$  multiplied by  $\varepsilon\nabla\Delta\operatorname{curl}\operatorname{u}$ to take the integral over  $\Omega\times (0,t)$ and $\Delta=\nabla \mathrm{div}-\mathrm{curl}\mathrm{curl}$, $\Delta \mathrm{curl}=-\mathrm{curl}\mathrm{curl}\mathrm{curl}$ to deduce that
\begin{align*}
\mu\int_{0}^{t}\int_{\Omega}(\rho-1)\vert\varepsilon^2\nabla\Delta \operatorname{curl}\mathbf u\vert^2dxds
=C_1+C_2+C_3+C_4+C_5+C_6,
\end{align*}
where
\begin{align*}
\begin{cases}
C_1=\int_{0}^{t}\int_{\Omega}\nabla(\varepsilon q \operatorname{curl}\operatorname{curl}\mathbf{u}_t)\cdot\varepsilon^2\nabla\Delta \operatorname{curl}\mathbf{u} dxds,\\
  C_2=\int_{0}^{t}\int_{\Omega}\nabla[\varepsilon q\nabla\times(\mathbf u\cdot\nabla \mathbf u)]\cdot\varepsilon^2\nabla\Delta \operatorname{curl}\mathbf{u} dxds,\\
  C_3=\int_{0}^{t}\int_{\Omega}\nabla\{\varepsilon q[\nabla\times(\mathbf{e}_3\times \mathbf u)]\}\cdot\varepsilon^2\nabla\Delta \operatorname{curl}\mathbf{u}dxds,\\
  C_4=\int_{0}^{t}\int_{\Omega}\nabla\{\varepsilon q[\nabla\times(\varepsilon q(\mathbf u_t+\mathbf u\cdot\nabla \mathbf u))]\}\cdot\varepsilon^2\nabla\Delta \operatorname{curl}\mathbf{u}dxds,\\
  C_5=\int_{0}^{t}\int_{\Omega}\nabla[\varepsilon q\nabla\times(\mathbf{e}_3\times q\mathbf u)]\cdot\varepsilon^2\nabla\Delta \operatorname{curl}\mathbf{u}dxds,\\
  C_6=\varepsilon\int_{0}^{t}\int_{\Omega}\nabla[\nabla\times(\rho\nabla\Delta\rho)]\cdot\varepsilon^2\nabla\Delta\operatorname{curl}
  \mathbf{u} dxds.
\end{cases}
\end{align*}
Therefore, we can combine $\Delta \mathrm{curl}=-\mathrm{curl}\mathrm{curl}\mathrm{curl}$ and Lemma \ref{lemma2.6}, Lemma \ref{lemma2.7}, to get the following estimate
\begin{align}\label{3.13}
\varepsilon\Vert\varepsilon \operatorname{curl}\mathbf u\Vert^2_{H^2}(t)+\varepsilon\Vert\varepsilon^2\operatorname{curl}\mathbf u\Vert^2_{L^2_t(H^3)}\leq C_0(M_0)\operatorname{exp}\{(t^\frac{1}{4}+\varepsilon^2)C_1(M(t))\}.
\end{align}
By combining $(\ref{3.10})$, $(\ref{3.12})$, $(\ref{3.13})$, we obtain estimate for  rotations of the velocity field.

By multiplying  $\varepsilon(\ref{1.3})_3$  by  $\varepsilon\Delta\theta$  to take the integral over  $\Omega\times (0,t)$  and combining the Lemma \ref{lemma2.6}, we get
\begin{align}\label{3.14}
\Vert(q,\theta)\Vert^2_{L^\infty_t(H^1)}+\Vert(\mathbf u,\varepsilon q,\varepsilon\theta)\Vert^2_{L^\infty_t(H^2)}+\Vert(\varepsilon q_t,\varepsilon\theta_t,\varepsilon \mathbf u_t)\Vert^2_{L^\infty_t(H^1)}+\Vert(\theta,\varepsilon \mathbf u_t,\varepsilon\theta_t)\Vert^2_{L^2_t(H^2)}\nonumber\\
\leq C_0(M_0)\operatorname{exp}\{(t^\frac{1}{8}+\varepsilon^{\frac{1}{2}})C_1(M(t))\}.
\end{align}
Next the equation  $(\ref{3.4})$  is derived with respect to $t$, which yields
\begin{align}\label{3.15}
\begin{cases}
q_{ttt}+\dfrac{1}{\varepsilon}\operatorname{div}\mathbf u_{tt}=\\
-(\mathbf u_{tt}\cdot\nabla q+2\mathbf u_t\cdot\nabla q_t+\mathbf u\cdot\nabla q_{tt}+q_{tt}\cdot \operatorname{div}\mathbf u+2q_t\operatorname{div}\mathbf u_t+q\operatorname{div}\mathbf u_{tt})~~~~~~~~~~~~~~~~~~~~~(t,x)\in(0,T]\times\Omega,\\
  \rho(\mathbf u_{ttt}+2\mathbf u_t\cdot\nabla \mathbf u_t+\mathbf u\cdot\nabla \mathbf u_{tt})+\dfrac{1}{\varepsilon}\nabla(q_{tt}+\theta_{tt})-\varepsilon\mu\Delta \mathbf u_{tt}-\varepsilon(\mu+\nu)\nabla \operatorname{div}\mathbf u_{tt}\\
  +\dfrac{1}{\varepsilon}(\mathbf{e}_3\times \mathbf u_{tt})=-2\rho_t(\mathbf u_{tt}+\mathbf u_t\cdot\nabla \mathbf u+\mathbf u\cdot\nabla \mathbf u_t)-\rho_{tt}(\mathbf u_t+\mathbf u\cdot\nabla \mathbf u_t)-\rho \mathbf u_{tt}\cdot\nabla \mathbf u\\
  -2\nabla(q_t\theta_t)-\nabla(q\theta_{tt})-\nabla (q_{tt}\theta)-\mathbf{e}_3\times(q\mathbf u)_{tt}+(\rho\nabla\Delta\rho)_{tt}~~~~~~~~~~~~~~~~~~~~~~~~~~~~~~~~~~~~~~~(t,x)\in(0,T]\times\Omega,\\
  (\rho \mathbf u)_{tt}\cdot\nabla\theta+2(\rho \mathbf u)_t\cdot\nabla\theta_t+\dfrac{1}{\varepsilon}\operatorname{div}\mathbf u_{tt}+\rho \mathbf u\cdot\nabla\theta_{tt}+\rho_{tt}\theta_t+2\rho_t\theta_{tt}+\rho\theta_{ttt}\\
  +\operatorname{div}\mathbf u_{tt}(\rho\theta+q)+2\operatorname{div}\mathbf u_t(\rho_t\theta+\rho\theta_t+q_t)+\operatorname{div}\mathbf u(\rho_{tt}\theta+2\rho_t\theta_t+\rho\theta_{tt}+q_{tt})\\
  =\lambda\Delta\theta_{tt}+(\varepsilon\Psi(\mathbf{u}^\varepsilon)+\Phi'(\rho^\varepsilon,\mathbf{u}^\varepsilon))_{tt}~~~~~~~~~~~~~~~~~~~~~~~~~~~~~~~~~~~~~~~~~~~~~~~~~~~~~~~~~~~~~~~~~~~~~~~~~(t,x)\in(0,T]\times\Omega,\\
  \mathbf u_{tt}\cdot{\mathbf n}=0,~~ {\mathbf n}\times \operatorname{curl}\mathbf u_{tt}|_=0~~~~~~~~~~~~~~~~~~~~~~~~~~~~~~~~~~~~~~~~~~~~~~~~~~~~~~~~~~~~~~~~~~~~~~~~~~~~~~~~~~~~(0,T]\times\partial\Omega,\\
  (q_{tt},\mathbf u_{tt})(0,x)=(q_{tt}(0),\mathbf u_{tt}(0))~~~~~~~~~~~~~~~~~~~~~~~~~~~~~~~~~~~~~~~~~~~~~~~~~~~~~~~~~~~~~~~~~~~~~~~~~~~~~~~~~x\in\Omega,
\end{cases}
\end{align}
where
\begin{align*}
q_{tt}(0):&=-\dfrac{1}{\varepsilon}\operatorname{div}\mathbf{u}_t(0)-q_t(0)\operatorname{div}\mathbf{u}_0-q_0\operatorname{div}\mathbf {u}_t(0)-\mathbf{u}_t(0)\cdot\nabla q_0-\mathbf{u}_0\cdot\nabla q_t(0),\\
\mathbf{u}_{tt}(0):&=[-\dfrac{1}{\varepsilon}(\nabla q_t(0)+\nabla\theta_t(0))+\varepsilon\mu\Delta\mathbf{u}_t(0)+\varepsilon(\mu+\nu)\nabla\operatorname{div}\mathbf{u}_t(0)-\rho_t(0)(\mathbf{u}_t(0)+\mathbf{u}_0\cdot\nabla\mathbf{u}_0)\\
&-\rho_0\mathbf{u}_t(0)\cdot\nabla \mathbf{u}_0-\nabla(q_t(0)\theta_0)-\nabla(q_0\theta_t(0))-\mathbf{e}_3\times(q_0\mathbf{u}_0)_t-\dfrac{1}{\varepsilon}(\mathbf{e}_3\times\mathbf{u}_t(0))
+(\rho_0\nabla\Delta\rho_0)_t]/\rho_0\\
&-\mathbf{u}_0\cdot\nabla\mathbf{u}_t(0),\\
\theta_{tt}(0):&=[\lambda\Delta\theta_t(0)+(\varepsilon\Psi(\mathbf{u}^\varepsilon)+\Phi'(\rho^\varepsilon,\mathbf{u}^\varepsilon))_t(0)-\frac{1}{\varepsilon}\operatorname{div}\mathbf{u}_t(0)-(\rho_0\mathbf{u}_0)_t\nabla\theta_0-\rho_0\mathbf{u}_0\cdot\nabla\theta_t(0)\\
&-\rho_t(0)\theta_t(0)-\operatorname{div}\mathbf u_t(0)(\rho_0\theta_0+q_0)-\operatorname{div}\mathbf{u}_0(\rho_t(0)\theta_0+q_t(0)+\rho_0\theta_t(0))]/\rho_0.
\end{align*}
Then we integrate the product of  $(\ref{3.15})_1$, $(\ref{3.15})_2$, $(\ref{3.15})_3$  with  $\varepsilon^4 q_{tt}$, $\varepsilon^4\mathbf{u}_{tt}$, $\varepsilon^4\theta_{tt}$ in $\Omega\times (0,t)$, and combine Lemma \ref{lemma2.6}, to obtain that
\begin{align}\label{3.16}
\Vert\varepsilon^2(q_{tt},\mathbf u_{tt},\theta_{tt})\Vert^2_{L^2}(t)+\varepsilon\Vert\varepsilon^2\mathbf u_{tt}\Vert^2_{L^2_t(H^1)}+\Vert\varepsilon^2\theta_{tt}\Vert^2_{L^2_t(H^1)}\leq C_0(M_0)\operatorname{exp}\{(t^\frac{1}{8}+\varepsilon)C_1(M(t))\}.
\end{align}
Next, we take the inner product of  $\varepsilon^2\partial_i\nabla(\ref{1.3})_2$, $\varepsilon^2\partial_i\nabla^2(\ref{1.3})_1$, $\varepsilon^2\partial_i\nabla^2(\ref{1.3})_3$ and $\varepsilon^2\partial_i\nabla^2\operatorname{div}\mathbf{u}$, $\varepsilon^2\partial_i\nabla^2 q$, $\varepsilon^2\partial_i\nabla^2\theta$ in $\Omega\times (0,t)$ to obtain that
\begin{align}\label{3.17}
\varepsilon\Vert\varepsilon^2\partial_i\nabla^2\operatorname{div}\mathbf u\Vert^2_{L^2_t(L^2)}+\Vert\varepsilon^2\partial_i\nabla^2q\Vert^2_{L^2}(t)+\Vert\varepsilon^2\partial_i\nabla^2
\theta\Vert^2_{L^2}(t)+\lambda\Vert\varepsilon^2\partial_i\nabla^3\theta\Vert^2_{L^2_t(L^2)}\nonumber\\
\leq C_0(M_0)\operatorname{exp}\{(t^\frac{1}{2}+\varepsilon)C_1(M(t))\}.
\end{align}
After that, we take the inner product of  $\partial_i\nabla(\ref{1.3})_1$, $\partial_i\nabla (\ref{1.3})_3$, $\partial_i(1.4)_2$ and $\varepsilon^2\partial_i\nabla q_t$, $\varepsilon^2\partial_i\nabla\theta_t$, $\varepsilon^2\partial_i\nabla\operatorname{div}\mathbf{u}$ in $\Omega\times (0,t)$ and combine $\Delta\mathbf{u}=\nabla\mathrm{div}\mathbf{u}-\mathrm{curlcurl}\mathbf{u}$, to obtain that
\begin{align}\label{3.18}
\Vert\varepsilon\partial_i\nabla q_t\Vert^2_{L^2_t(L^2)}+\Vert\varepsilon\partial_i\nabla \theta_t\Vert^2_{L^2_t(L^2)}+\varepsilon\Vert\varepsilon\partial_i\nabla \operatorname{div}\mathbf u\Vert^2_{L^2_t(L^2)}\leq C_0(M_0)\operatorname{exp}\{(t^\frac{1}{8}+\varepsilon^{\frac{1}{2}})C_1(M(t))\}.
\end{align}
Combining  $(\ref{3.17})$  and  $(\ref{3.18})$  yields
\begin{align}\label{3.19}
&\varepsilon\Vert\varepsilon^2\nabla^3\operatorname{div}\mathbf u\Vert^2_{L^2_{t}(L^2)}+\Vert(\varepsilon\nabla^2q_t,\varepsilon\nabla^2\theta_t,\varepsilon\nabla^3q,\varepsilon\nabla^3\theta,\varepsilon^2\nabla^4\theta)\Vert^2_{L^2_t(L^2)}+\Vert\varepsilon \mathbf u\Vert^2_{L^2_{t}(H^3)}+\varepsilon\Vert\varepsilon\nabla^2\operatorname{div}\mathbf u\Vert^2_{L^2}(t)\nonumber\\
&+\Vert(\varepsilon^2\nabla^3q,\varepsilon^2\nabla^3\theta)\Vert^2_{L^2}(t)\leq C_0(M_0)\operatorname{exp}\{(t^\frac{1}{8}+\varepsilon)C_1(M(t))\}.
\end{align}
Combining  $(\ref{3.14})$, $(\ref{3.16})$, $(\ref{3.19})$ gives higher order derivative estimates.

Next, $\varepsilon\nabla(1.3)_2$ integrates over $\Omega\times{(0,t)}$, we obtain
\begin{align*}
\Vert q\Vert^2_{L^2_t(H^2)}\leq C_0(M_0)\mathrm{exp}\{(t^\frac{1}{8}+\varepsilon^{\frac{1}{2}})C_1(M(t))\}.
\end{align*}

Combining all the above results gives
\begin{align*}
M^\varepsilon(t)\leq C_0(M^\varepsilon_0)\mathrm{exp}\{(t^\frac{1}{8}+\varepsilon^{\frac{1}{2}})C_1(M^\varepsilon(t))\},
\end{align*}
where $C_1(\cdot)$ and $C_0(\cdot)$ are both positive continuous decreasing functions.

There is an a priori estimate that when $\varepsilon\to 0$
\begin{align}
q^\varepsilon\stackrel{w}{\longrightarrow} q~~~&\operatorname{in} L^2(0,T;H^2(\Omega));\qquad q^\varepsilon\stackrel{w*}{\longrightarrow} q~~~\operatorname{in} L^\infty(0,T;H^1(\Omega)),\\
\mathbf u^\varepsilon\stackrel{w}{\longrightarrow} \mathbf u~~~&\operatorname{in} L^2(0,T;H^3(\Omega));\qquad \mathbf u^\varepsilon\stackrel{w*}{\longrightarrow} \mathbf u~~~\operatorname{in} L^\infty(0,T;H^2(\Omega)),\\
\theta^\varepsilon\stackrel{w}{\longrightarrow} \theta~~~&\operatorname{in}
L^2(0,T;H^3(\Omega));\qquad \theta^\varepsilon\stackrel{w*}{\longrightarrow} \theta~~~\operatorname{in} L^\infty(0,T;H^1(\Omega)),
\end{align}
and 
\begin{align}
&\rho^\varepsilon=1+\varepsilon q^\varepsilon\to 1~~~\operatorname{in} L^2(0,T;H^2(\Omega))\cap L^\infty(0,T;H^1(\Omega)),\\
&\Theta^\varepsilon=1+\varepsilon\theta^\varepsilon\to 1~~~\operatorname{in} L^2(0,T;H^2(\Omega))\cap L^\infty(0,T;H^1(\Omega)).
\end{align}

\subsection{Acoustic analysis}\label{s3.2}

When  $\alpha=1$, combine $(\ref{1.3})_2$ and   $(\ref{1.3})_3$ to form the sonic equation
\begin{align}\label{3.25}
\begin{cases}
\varepsilon\partial_t\sigma^\varepsilon+\operatorname{div}\mathbf U^\varepsilon=\varepsilon f^\varepsilon,\\
 \varepsilon \mathbf U^\varepsilon_t+(\mathbf{e}_3\times \mathbf U^\varepsilon+\nabla_x\sigma^\varepsilon)=\varepsilon l^\varepsilon,
\end{cases}
\end{align}
where
\begin{align*}
\mathbf U^\varepsilon&=\rho^\varepsilon \mathbf u^\varepsilon,\\
\sigma^\varepsilon&=q^\varepsilon+\theta^\varepsilon,~~~~~~~~~~~~~~~~\sigma^\varepsilon(0,\cdot)=\sigma^\varepsilon_0=q^\varepsilon_0+\theta^\varepsilon_0,\\
f^\varepsilon&=\lambda\Delta\theta^\varepsilon+\mathbf u^\varepsilon \operatorname{div}q^\varepsilon+q_{t}^\varepsilon-\varepsilon\theta^\varepsilon_{t} q^\varepsilon_{t}-\rho^\varepsilon \mathbf u^\varepsilon\cdot\nabla\theta^\varepsilon-\theta^\varepsilon \operatorname{div}\mathbf u^\varepsilon-\varepsilon q^\varepsilon\theta^\varepsilon \operatorname{div}\mathbf u^\varepsilon+\varepsilon(\nu(\operatorname{div}\mathbf{u}^\varepsilon)^2\\
&~~~~+2\mu D(\mathbf{u}^\varepsilon):\nabla\mathbf{u}^\varepsilon)+\rho^\varepsilon\Delta q^\varepsilon\operatorname{div}\mathbf{u}^\varepsilon+\frac{\varepsilon}{2}\vert\nabla q^\varepsilon\vert^2\operatorname{div}\mathbf{u}^\varepsilon-\varepsilon(\nabla q^\varepsilon\otimes\nabla q^\varepsilon):\nabla\mathbf{u}^\varepsilon,\\
l^\varepsilon&=-\operatorname{div}(\rho^\varepsilon \mathbf u^\varepsilon\otimes \mathbf u^\varepsilon)-\nabla(q^\varepsilon\theta^\varepsilon)+\varepsilon(\mu+\nu)\nabla\operatorname{div}\mathbf{u}^\varepsilon+\varepsilon\mu\Delta\mathbf{u}^\varepsilon+\rho^\varepsilon\nabla\Delta\rho^\varepsilon.
\end{align*}
For any $\varphi\in([0,T)\times\overline{\Omega})$, this follows from $(\ref{3.25})_1$ and $(\ref{3.25})_2$
\begin{align*}
&\int_{0}^{T}\int_{\Omega}\varepsilon \sigma^\varepsilon\cdot\partial_t\varphi+\mathbf{U}^\varepsilon\cdot \nabla_x\varphi dxdt=-\int_{\Omega}\varepsilon \sigma^\varepsilon_0\cdot\varphi(0,\cdot)dx+\int_{0}^{T}\varepsilon\left<f^\varepsilon,\varphi\right>dt,\\
&\int_{0}^{T}\int_{\Omega}-\varepsilon\rho^\varepsilon \mathbf u^\varepsilon\cdot\partial_t\varphi+\mathbf{e}_3\times \mathbf{U}^\varepsilon\cdot\varphi
-\sigma^\varepsilon\operatorname{div}\varphi dxdt=\int_{\Omega}\varepsilon\rho_0^\varepsilon \mathbf u_0^\varepsilon\cdot\varphi(0,\cdot)dx+\int_{0}^{T}\varepsilon\left<l^\varepsilon,\varphi\right>dt,
\end{align*}
where
\begin{align*}
&\left<f^\varepsilon,\varphi\right>\\
&=\int_\Omega\lambda\nabla\theta^\varepsilon\operatorname{div}\varphi+\mathbf{u}^\varepsilon\operatorname{div}q^\varepsilon\varphi
+q_t^\varepsilon\varphi-\varepsilon\theta_t^\varepsilon q_t^\varepsilon\varphi-\rho^\varepsilon\mathbf{u}^\varepsilon\cdot\nabla\theta^\varepsilon\cdot\varphi-\theta^\varepsilon\operatorname{div}
\mathbf{u}^\varepsilon\varphi-\varepsilon q^\varepsilon\theta^\varepsilon\operatorname{div}\mathbf{u}^\varepsilon\varphi\\
&~~~~+[\varepsilon(\nu(\operatorname{div}\mathbf{u}^\varepsilon)^2+2\mu D(\mathbf{u}^\varepsilon):\nabla\mathbf{u}^\varepsilon)+\rho^\varepsilon\Delta q^\varepsilon\operatorname{div}\mathbf{u}^\varepsilon+\frac{\varepsilon}{2}\vert\nabla q^\varepsilon\vert^2\operatorname{div}\mathbf{u}^\varepsilon-\varepsilon(\nabla q^\varepsilon\otimes\nabla q^\varepsilon):\nabla\mathbf{u}^\varepsilon]\cdot\varphi dx\\
&=\int_\Omega(f_1^\varepsilon\operatorname{div}\varphi+f_2^\varepsilon\cdot\varphi)dx.\\
&\left<l^\varepsilon,\varphi\right>\\
&=\int_{\Omega}(\rho^\varepsilon \mathbf u^\varepsilon\otimes \mathbf u^\varepsilon):\nabla_x\varphi+(q^\varepsilon\theta^\varepsilon) \operatorname{div}\varphi-2\mu D\mathbf{u}^\varepsilon:\nabla\varphi+\nu(\operatorname{div}\mathbf{u}^\varepsilon)^2\cdot\varphi+\rho^\varepsilon\nabla\Delta\rho^\varepsilon\cdot\varphi dx\\
&=\int_\Omega(l_1^\varepsilon:\nabla\varphi+l_2^\varepsilon\operatorname{div}\varphi+l_3^\varepsilon\cdot\varphi)dx.
\end{align*}
By the Sobolev embedding theorem, one obtains  $f_1^\varepsilon$  is  bounded in  $L^\infty(0,T;L^2)\cap L^2(0,T;H^1)$,~$f_2^\varepsilon$  is bounded in  $L^\infty(0,T;H^2)\cap L^2(0,T;L^1)$;~$l_1^\varepsilon$  is bounded in  $L^\infty(0,T;L^1)\cap L^2(0,T;L^1)$,~$l_2^\varepsilon$  is bounded in  $L^\infty(0,T;L^1)\cap L^2(0,T;L^1)$,~$l_3^\varepsilon$  is bounded in  $L^2(0,T;L^2)$.

Consider an operator  $\mathcal{B}$  defined, formally, in  $L^2(\Omega)\times L^2(\Omega;R^3)$,
\begin{align*}
\mathcal{B}:\begin{bmatrix}
\sigma^\varepsilon\\
\mathbf U^\varepsilon
\end{bmatrix}
\longmapsto
\begin{bmatrix}
\operatorname{div}\mathbf U^\varepsilon\\
\mathbf{e}_3\times \mathbf U^\varepsilon+\nabla_x\sigma^\varepsilon
\end{bmatrix}.
\end{align*}
We study the point spectrum of  $\mathcal{B}$, we look for solutions to the eigenvalue problem
\begin{align}
\mathcal{B}\begin{bmatrix}
\sigma^\varepsilon\\
\mathbf U^\varepsilon
\end{bmatrix}
=
\begin{bmatrix}
\operatorname{div}\mathbf U^\varepsilon\\
e_3\times\mathbf U^\varepsilon+\nabla_x\sigma^\varepsilon
\end{bmatrix}
=\lambda
\begin{bmatrix}
\sigma^\varepsilon\\
\mathbf U^\varepsilon
\end{bmatrix},
\end{align}
or in terms of Fourier variables
\begin{align*}
i(\sum_{j=1}^2\xi_j\hat{\mathbf{U}^\varepsilon_j}+k\hat{\mathbf{U}^\varepsilon_3})-\lambda\hat{\sigma^\varepsilon}=0,~~~~i[\xi_1,\xi_2,k](\hat{\sigma^\varepsilon})-[\hat{\mathbf{U}^\varepsilon_2},
-\hat{\mathbf{U}^\varepsilon_1},0]-\lambda\hat{\mathbf{U}^\varepsilon}=0.
\end{align*}
We can find that
\begin{align*}
\lambda=-\mu,~~~\mu=\frac{1+\vert\xi\vert^2+k^2\pm\sqrt{(1+\vert\xi\vert^2+k^2)^2-4k^2}}{2}.
\end{align*}
We can obtain that when  $k=0$, the only eigenvalue is  $\lambda=0$. Therefore, the eigenvector space coincides with the kernel space of  $\mathcal{B}$.
\begin{align*}
&\operatorname{ker}(\mathcal{B})=\{[\sigma^\varepsilon,\mathbf U^\varepsilon]| \sigma^\varepsilon=\sigma^\varepsilon(x_1,x_2),\nonumber\\
\mathbf U^\varepsilon=[\mathbf{U}^\varepsilon_1(x_1,x_2),~\mathbf{U}^\varepsilon_2&(x_1,x_2),~\mathbf{U}^\varepsilon_3(x_1,x_2)],~\operatorname{div}_h\mathbf{U}^\varepsilon_h=0,~\nabla_h\sigma^\varepsilon=
[\mathbf{U}^\varepsilon_2,-\mathbf{U}^\varepsilon_1]\}.
\end{align*}

First, the Hilbert space is defined
\begin{align*}
H=H_M=\{[\sigma,\mathbf U]|\hat{\sigma}(\xi_h,k)=0,~\hat{\mathbf U}(\xi_h,k)=0, \operatorname{where} \vert\xi_h\vert+\vert k\vert\textgreater M\}.
\end{align*}
And let the orthogonal projection on  $H_M$  be given by
\begin{align*}
P_M:L^2(\Omega)\times L^2(\Omega;R^3)\to H_M,
\end{align*}

Combined with  $(\ref{3.25})$, this yields
\begin{align}\label{3.27}
\varepsilon\frac{d}{dt}\begin{bmatrix}
\sigma^\varepsilon_M\\
\mathbf{U}^\varepsilon_M
\end{bmatrix}+\mathcal{B}
\begin{bmatrix}
\sigma^\varepsilon_M\\
\mathbf{U}^\varepsilon_M\end{bmatrix}
=\varepsilon\begin{bmatrix}
f^\varepsilon_M\\
l^\varepsilon_M
\end{bmatrix},
\end{align}
where  $\mathcal{A}=i\mathcal{B}$, $(f^\varepsilon_M,l^\varepsilon_M)\in H_M^{*}\cong H_M$, when  $(s,w)\in H_M$, it follows that
\begin{align*}
\left<\begin{bmatrix}
f^\varepsilon_M\\
l^\varepsilon_M
\end{bmatrix}
,
\begin{bmatrix}
s\\
w
\end{bmatrix}\right>
=-\int_{\Omega}\begin{bmatrix}
f_1^\varepsilon\operatorname{div}s+f_2^\varepsilon\cdot s\\
l_1^\varepsilon:\nabla w+l_2^\varepsilon\operatorname{div}w+l_3^\varepsilon\cdot w
\end{bmatrix}dx.
\end{align*}
Combined with the a priori estimate, as $\varepsilon\to0$
\begin{align}\label{3.28}
\bigg\Vert
\begin{bmatrix}
f^\varepsilon\\
l^\varepsilon
\end{bmatrix}\displaystyle\bigg\Vert_{L^2(0,T;P)}\leq C(M).
\end{align}

We defined
\begin{align*}
Q:L^2(\Omega)\times L^2(\Omega;R^3)\to \operatorname{ker}\mathcal B,
\end{align*}
the orthogonal projection onto the null space of  $\mathcal{B}$, from $[P_M,Q]=0$, we can obtain
\begin{align*}
\bigg\Vert\sqrt{C}Q^{\perp}
\begin{bmatrix}
\sigma^\varepsilon_M\\
\mathbf{U}^\varepsilon_M
\end{bmatrix}
\bigg\Vert^2_{H_M}&=\bigg<CQ^{\perp}
\begin{bmatrix}
\sigma^\varepsilon_M\\
\mathbf{U}^\varepsilon_M
\end{bmatrix},
\begin{bmatrix}
\sigma^\varepsilon_M\\
\mathbf{U}^\varepsilon_M
\end{bmatrix}\bigg>_{H_M}=\int_{\Omega}
\chi\vert Q^{\perp}(\sigma^\varepsilon_M,\mathbf{U}^\varepsilon_M)\vert^2 dx.
\end{align*}
Combined  with the  Rage  theorem, this yields that when  $\varepsilon\to0$,
\begin{align}\label{3.29}
Q^{\perp}(\sigma^\varepsilon_M,\mathbf{U}^\varepsilon_M)\to0,~~~\operatorname{in}~L^2([0,T]\times K; R^4),
\end{align}
for any compact  $K\subset\overline{\Omega}$ and any fixed  $M$.

By the embedding theorem and  Ascoli-Arzela  theorem, we obtain that when $\varepsilon\to0$,
\begin{align}\label{3.30}
Q(\sigma^\varepsilon_M,\mathbf{U}^\varepsilon_M)\to(\sigma_M,\overline{\rho}\mathbf{u}_M),~~~\operatorname{in}~L^2([0,T]\times K; R^4).
\end{align}

By virtue of  $(\ref{3.29})$ and (\ref{3.30}), we obtain
\begin{align*}
P_M[\mathbf{u}^\varepsilon]\to P_M[\mathbf u],~~~\operatorname{in}~L^2([0,T]\times K;R^3),
\end{align*}
for any fixed $M$,  and compactness of the embedding $H^{1}(M) \hookrightarrow\hookrightarrow L^{2}(M)$, we get
\begin{align}\label{3.31}
\mathbf{u}^\varepsilon\to \mathbf u,~~~\operatorname{in}~~L^2([0,T]\times K;R^3),~~~\operatorname{for~any~compact}~~ K\subset\Omega.
\end{align}
\subsection{Identifying the limit system}

Let us set $\varphi=[\nabla^\perp_h\psi,0],~~\psi\in C^\infty_c([0,T)\times\Omega)$ from $(\ref{1.3})_2$, we get
\begin{align}\label{3.32}
&\int_{0}^{T}\int_{\Omega}\rho^\varepsilon \mathbf u^\varepsilon\cdot\partial_t\varphi+\rho^\varepsilon \mathbf u^\varepsilon\otimes \mathbf u^\varepsilon:\nabla_x\varphi-\frac{1}{\varepsilon}\rho^\varepsilon[\mathbf u^\varepsilon]_h\cdot\nabla_x\psi dxdt\nonumber\\
&=-\int_{\Omega}\rho_0^\varepsilon \mathbf u_0^\varepsilon\cdot\varphi(0,\cdot)dx+\int_{0}^{T}\int_{\Omega}\varepsilon(\mu+\nu)\operatorname{div}\mathbf{u}^\varepsilon:
\nabla_x\varphi-\varepsilon\mu\Delta\mathbf{u}^\varepsilon\cdot\varphi-\varepsilon\rho^\varepsilon\nabla\Delta q^\varepsilon
\cdot\varphi dxdt.
\end{align}
From  $(\ref{1.3})_3$, we obtain
\begin{align}\label{3.33}
&\int_{0}^{T}\int_{\Omega}\sigma^\varepsilon\cdot\partial_t\psi+\frac{1}{\varepsilon}\rho^\varepsilon[\mathbf u^\varepsilon]_h\cdot \nabla_x\psi dxdt\nonumber\\
&=-\int_{\Omega}
\sigma^\varepsilon_0\psi(0,\cdot)dx+\int_{0}^{T}\int_{\Omega}\lambda\nabla\theta^\varepsilon\cdot\operatorname{div}\psi dxdt+\int_{\Omega}\int_0^t[\varepsilon(\nu(\operatorname{div}\mathbf{u}^\varepsilon)^2\nonumber+2\mu D(\mathbf{u}^\varepsilon):\nabla\mathbf{u}^\varepsilon)\\
&~~~~+\mathbf{u}^\varepsilon\operatorname{div}q^\varepsilon+q^\varepsilon_t-\varepsilon q^\varepsilon_t\theta^\varepsilon_t-\rho^\varepsilon\mathbf{u}^\varepsilon\nabla\theta^\varepsilon-\theta^\varepsilon\operatorname{div}\mathbf{u}^\varepsilon\nonumber-\varepsilon q^\varepsilon\theta^\varepsilon\operatorname{div}\mathbf{u}^\varepsilon+\rho^\varepsilon\Delta q^\varepsilon\operatorname{div}\mathbf{u}^\varepsilon+\frac{\varepsilon}{2}\vert\nabla q^\varepsilon\vert^2\operatorname{div}\mathbf{u}^\varepsilon\nonumber\\
&~~~~-\varepsilon(\nabla q^\varepsilon\otimes\nabla q^\varepsilon):\nabla\mathbf{u}^\varepsilon]\cdot\psi dtdx.
\end{align}
Combining  $(\ref{3.32})$ and  $(\ref{3.33})$  yields that  when  $\varepsilon\to0$
\begin{align}\label{3.34}
&\int_{0}^{T}\int_{\Omega}(\bar{\rho}\mathbf u_h\cdot\partial_t\nabla_h^\perp\psi+\bar{\rho}[\mathbf u_h\otimes \mathbf u_h]:\nabla_h(\nabla_h^\perp\psi)+\sigma\cdot\partial_t\psi)dxdt\nonumber\\
&=-\int_{\Omega}(\bar{\rho}(\int_{0}^{1}\mathbf u_{0,h}dx_3)\cdot\nabla^\perp_h\psi(0,\cdot)+(\int_{0}^{1}\sigma_0dx_3)\psi(0,\cdot))dx_h.
\end{align}

Combining the diagnostic equations,  $\mathbf{u}_h=\nabla^\perp_h\sigma$ and  $(\ref{3.34})$ yields the convergence result of Theorem \ref{th1.2}.
\section{Limit problem of capillary effect when $\alpha=0$}\label{s4}

This section proves the Theorem \ref{th1.3} of this thesis. 

When $\alpha=0$, the capillary term is $\dfrac{1}{\varepsilon^{2(1-\alpha)}}\rho^\varepsilon\nabla\Delta\rho^\varepsilon=\dfrac{1}{\varepsilon^2}\rho^\varepsilon\nabla\Delta\rho^\varepsilon$.
We let $\mathbf w$ and $\nabla\pi$ denote $\mathbf{u}^\varepsilon$ and $\frac{1}{\varepsilon}[(1+\varepsilon q^\varepsilon)\nabla\theta^\varepsilon+(1+\varepsilon\theta^\varepsilon)\nabla q^\varepsilon-(1+\varepsilon q^\varepsilon)\nabla\Delta q^\varepsilon-(1+\varepsilon q^\varepsilon)\nabla((I-\Delta)q^\varepsilon+\theta^\varepsilon)]$ in the limit, $(\mathbf w,\pi)$ is a solution to the incompressible Navier-Stokes-Korteweg equations $(\ref{4.1})$
\begin{align}\label{4.1}
\begin{cases}
\mathbf{w}_t+\mathbf{w}\cdot\nabla_h\mathbf{w}+\nabla_h\pi=0,\\
 \operatorname{div}\mathbf{w}=0,\\
 \mathbf{w}(0,x)=\mathbf{w}_0,~~~~~\operatorname{div}\mathbf{w}_0=0.
\end{cases}
\end{align}
Construct an approximate solution $(q_\varepsilon,\mathbf{u}_\varepsilon,\theta_\varepsilon):=(\varepsilon^2\pi,\mathbf w,\varepsilon\pi)$ satisfies
\begin{align}\label{4.2}
\begin{cases}
q_\varepsilon+\mathbf{u}_\varepsilon\cdot\nabla q_\varepsilon+\frac{1}{\varepsilon}(1+\varepsilon q_\varepsilon)\operatorname{div}\mathbf{u}_\varepsilon=R_1,\\
 (1+\varepsilon q_\varepsilon)(\mathbf{u}_\varepsilon+\mathbf{u}_\varepsilon\cdot\nabla\mathbf{u}_\varepsilon)+\frac{1}{\varepsilon}[(1+\varepsilon q_\varepsilon)\nabla\theta_\varepsilon+(1+\varepsilon\theta_\varepsilon)\nabla q_\varepsilon]+\frac{1}{\varepsilon}(1+\varepsilon q_\varepsilon)(e_3\times\mathbf{u}_\varepsilon)=\varepsilon\mu\Delta\mathbf{u}_\varepsilon+R_2,\\
 (1+\varepsilon q_\varepsilon)(\theta_\varepsilon+\mathbf{u}_\varepsilon\cdot\nabla\theta_\varepsilon)+\frac{1}{\varepsilon}(1+\varepsilon q_\varepsilon)(1+\varepsilon \theta_\varepsilon)\operatorname{div}\mathbf{u}_\varepsilon=\lambda\Delta\theta_\varepsilon+R_3,
\end{cases}
\end{align}
where
\begin{align*}
&R_1=\varepsilon^2(\pi_t+\mathbf{w}\cdot\nabla\pi),~~~R_2=\varepsilon(1+\varepsilon^3\pi)\nabla\Delta\pi-(1+\varepsilon^4\pi)\nabla\pi-\varepsilon\mu\Delta\pi,\\
&R_3=\varepsilon(1+\varepsilon^3\pi)(\pi_t+\mathbf{w}\cdot\nabla\pi)-\varepsilon\lambda\Delta\pi.
\end{align*}
For any $t\in[0,t^\ast]$,
\begin{align}\label{4.3}
\Vert R_1(t)\Vert_{H^s(\Omega)}\leq C\varepsilon^2,~~~\Vert R_2(t)\Vert_{H^s(\Omega)}\leq C\varepsilon+C,~~~\Vert R_3(t)\Vert_{H^s(\Omega)}\leq C\varepsilon.
\end{align}
\subsection{Error estimate}

In this section, we prove Theorem \ref{th1.3} via the energy error estimation method. First, we introduce the error variables 
\begin{align*}(q, \mathbf{u}, \theta)=(q^\varepsilon-q_\varepsilon, \mathbf{u}^\varepsilon-\mathbf{u}_\varepsilon, \theta^\varepsilon-\theta_\varepsilon),~t\in[0,\operatorname{min}\{T^{\ast}, T_\varepsilon\}].
\end{align*}
Then the difference between system $(\ref{1.3})$ and $(\ref{4.2})$ reads as follows:
\begin{align}\label{4.4}
\begin{cases}
q_t+\mathbf{u}^\varepsilon\cdot\nabla q+\frac{1}{\varepsilon}(1+\varepsilon q^\varepsilon)\operatorname{div}\mathbf{u}=-\mathbf{u}\cdot\nabla q_\varepsilon-R_1,\\
 \mathbf{u}_t+\mathbf{u}\cdot\nabla\mathbf{u}_\varepsilon+\mathbf{u}^\varepsilon\cdot\nabla\mathbf{u}+
 \frac{1}{\varepsilon}(\nabla\theta+\frac{1+\varepsilon\theta^\varepsilon}{1+\varepsilon q^\varepsilon}\nabla q)+\frac{1}{\varepsilon}(e_3\times\mathbf{u})=-\frac{1}{\varepsilon}(\frac{1+\varepsilon\theta^\varepsilon}{1+\varepsilon q^\varepsilon}-\frac{1+\varepsilon\theta_\varepsilon}{1+\varepsilon q_\varepsilon})\nabla q_\varepsilon\\
 +\frac{\varepsilon\mu\Delta\mathbf{u}}{1+\varepsilon q^\varepsilon}+\frac{\varepsilon(\mu+\nu)}{1+\varepsilon q^\varepsilon}\nabla\operatorname{div}\mathbf{u}+\frac{1}{\varepsilon}\nabla\Delta q^\varepsilon-\frac{R_2}{1+\varepsilon q_\varepsilon}+(\frac{\varepsilon\mu}{1+\varepsilon q^\varepsilon}-\frac{\varepsilon\mu}{1+\varepsilon q_\varepsilon})\Delta\mathbf{u},\\
 \theta_t+\mathbf{u}^\varepsilon\cdot\nabla\theta+\frac{1}{\varepsilon}(1+\varepsilon\theta^\varepsilon)\operatorname{div}\mathbf{u}+\mathbf u\cdot\nabla\theta_\varepsilon=\frac{Q}{1+\varepsilon q^\varepsilon}+(\frac{\lambda}{1+\varepsilon q^\varepsilon}-\frac{\lambda}{1+\varepsilon q_\varepsilon})\Delta\theta_\varepsilon+\frac{\lambda\Delta\theta}{1+\varepsilon q^\varepsilon}-\frac{R_3}{1+\varepsilon q_\varepsilon},
\end{cases}
\end{align}
where $Q=\varepsilon\Psi(\mathbf{u}^\varepsilon)+\kappa((1+\varepsilon q^\varepsilon)\Delta q^\varepsilon+\frac{\varepsilon|\nabla q^\varepsilon|^2}{2})\mathbf{u}^\varepsilon-\kappa\varepsilon(\nabla q^\varepsilon\otimes\nabla q^\varepsilon):\nabla\mathbf{u}^\varepsilon$, taking the operator $D^\alpha$ for equation $(\ref{4.4})$, we obtain
\begin{align}\label{4.5}
\partial_t q^\alpha+\mathbf{u}^\varepsilon\cdot\nabla q^\alpha+\frac{1}{\varepsilon}(1+\varepsilon q^\varepsilon)\operatorname{div}\mathbf{u}^\alpha+[D^\alpha,\mathbf{u}^\varepsilon]\cdot\nabla q^\alpha=-[D^\alpha, \frac{1}{\varepsilon}(1+\varepsilon q^\varepsilon)]\operatorname{div}\mathbf{u}-D^\alpha(\mathbf{u}\cdot\nabla q_{\varepsilon})-D^\alpha R_1,
\end{align}
\begin{align}\label{4.6}
&\partial_t\mathbf{u}^\alpha+\mathbf{u}^\varepsilon\cdot\nabla\mathbf{u}^\alpha+\frac{1}{\varepsilon}(\nabla\theta^\alpha+
\frac{1+\varepsilon\theta^\varepsilon}{1+\varepsilon q^\varepsilon}\nabla q^\alpha)+\frac{1}{\varepsilon}e_3\times\mathbf u^\alpha-\frac{\varepsilon\mu}{1+\varepsilon q^\varepsilon}\Delta\mathbf{u}^\alpha\nonumber\\
&=\frac{\varepsilon(\mu+\nu)}{1+\varepsilon q^\varepsilon}\nabla\operatorname{div}\mathbf u^\alpha-D^\alpha(\mathbf{u}\cdot\nabla\mathbf{u}_\varepsilon)-\frac{1}{\varepsilon}D^\alpha((\frac{1+\varepsilon\theta^\varepsilon}{1+\varepsilon
q^\varepsilon}-\frac{1+\varepsilon\theta_\varepsilon}{1+\varepsilon q_\varepsilon})\nabla q_\varepsilon)+[D^\alpha,\frac{\varepsilon\mu}{1+\varepsilon q^\varepsilon}]\Delta\mathbf{u}\nonumber\\
&~~~~+[D^\alpha,\frac{\varepsilon(\mu+\nu)}{1+\varepsilon q^\varepsilon}]\nabla\operatorname{div}\mathbf u+D^\alpha((\frac{\varepsilon\mu}{1+\varepsilon q^\varepsilon}-\frac{\varepsilon\mu}{1+\varepsilon q_\varepsilon})\Delta\mathbf{u}_\varepsilon)+\frac{1}{\varepsilon}D^\alpha\nabla\Delta q^\varepsilon-D^\alpha(\frac{R_2}{1+\varepsilon q_\varepsilon})\nonumber\\
&~~~~-[D^\alpha,\frac{1}{\varepsilon}\frac{1+\varepsilon\theta^\varepsilon}{1+\varepsilon q^\varepsilon}]\nabla q-[D^\alpha,\mathbf{u}^\varepsilon]\cdot\nabla\mathbf u,
\end{align}
\begin{align}\label{4.7}
&\partial_t\theta^\alpha+\mathbf{u}^\varepsilon\cdot\nabla\theta^\alpha+\frac{1}{\varepsilon}(1+\varepsilon\theta^\varepsilon)\operatorname{div}
\mathbf{u}^\alpha-\frac{\lambda}{1+\varepsilon q^\varepsilon}\Delta\theta^\alpha\nonumber\\
&=-[D^\alpha,\mathbf{u}^\varepsilon]\cdot\nabla\theta-[D^\alpha,\frac{1}{\varepsilon}(1+\varepsilon\theta^\varepsilon)]\operatorname{div}\mathbf u-D^\alpha(\mathbf{u}\cdot\nabla\theta_\varepsilon)+D^\alpha(\frac{Q}{1+\varepsilon q^\varepsilon})+[D^\alpha,\frac{\lambda}{1+\varepsilon q^\varepsilon}]\Delta\theta\nonumber\\
&~~~~+D^\alpha((\frac{\lambda}{1+\varepsilon q^\varepsilon}-\frac{\lambda}{1+\varepsilon q_\varepsilon})\Delta\theta_\varepsilon)-D^\alpha(\frac{R_3}{1+\varepsilon q_\varepsilon}).
\end{align}
Let $E=(q,\mathbf u,\theta)$, then $\Vert E\Vert_s=\Vert(q,\mathbf u,\theta)\Vert_s$.
By combining lemma \ref{lemma2.2} and Sobolev's inequality, we obtain
\begin{align}\label{4.8}
\Vert\partial_x^\beta q^\varepsilon\Vert_{L^\infty}\leq C\Vert E\Vert_{H^s(\Omega)}+C\varepsilon^2,~~~\Vert\partial_x^\beta\boldsymbol{u}^\varepsilon\Vert_{L^\infty}\leq C\Vert E\Vert_{H^s(\Omega)}+C,\nonumber\\
\Vert\partial_x^\beta \theta^\varepsilon\Vert_{L^\infty}\leq C\Vert E\Vert_{H^s(\Omega)}+C\varepsilon,~~~\Vert\partial_x^\gamma q^\varepsilon\Vert_{L^\infty}\leq C\Vert E\Vert_{H^s(\Omega)}+C\Vert\nabla q\Vert_{H^s(\Omega)}+C\varepsilon^2,
\end{align}
Combining $(\ref{1.3})_1$, $(\ref{1.3})_3$, we obtain 
\begin{align}\label{4.9}
\Vert\varepsilon q_t^\varepsilon\vert_{L^\infty}\leq C\Vert E\Vert_{H^s(\Omega)}+C\Vert E\Vert_{H^s(\Omega)}^2+C\varepsilon^2+C,\nonumber\\
\Vert\varepsilon \theta_t^\varepsilon\vert_{L^\infty}\leq C\Vert E\Vert_{H^s(\Omega)}+C\Vert E\Vert_{H^s(\Omega)}^2+C\Vert E\Vert_{H^s(\Omega)}^3+C\varepsilon+C.
\end{align}

Taking $(\ref{4.5})$, $(\ref{4.6})$, $(\ref{4.7})$ and $\frac{1+\varepsilon\theta^\varepsilon}{1+\varepsilon q^\varepsilon}q^\alpha$, $(1+\varepsilon q^\varepsilon)\mathbf{u}^\alpha$, $\frac{1+\varepsilon q^\varepsilon}{1+\varepsilon\theta^\varepsilon}\theta^\alpha$ using the partition integral to add the three equations.
By adding these three equations, and combining them with lemma \ref{lemma2.2}, lemma \ref{lemma2.5}, (\ref{4.8}), (\ref{4.9}), Cauchy's inequality and Young's inequality, we get
\begin{align}\label{4.10}
&\frac{1}{2}\frac{d}{dt}\sum\limits_{|\alpha|\leq s}\int_\Omega(\frac{1+\varepsilon\theta^\varepsilon}{1+\varepsilon q^\varepsilon}\vert q^\alpha\vert^2+(1+\varepsilon q^\varepsilon)\vert\mathbf{u}^\alpha\vert^2+\frac{1+\varepsilon q^\varepsilon}{1+\varepsilon\theta^\varepsilon}\vert\theta^\alpha\vert^2+\vert\nabla q^\alpha\vert^2)dx\nonumber\\
&+\frac{\varepsilon\mu}{2}\Vert\nabla\mathbf{u}\Vert_{H^s(\Omega)}^2+\frac{\varepsilon(\mu+\nu)}{2}\Vert\operatorname{div}\mathbf u\Vert_{H^s(\Omega)}^2+\frac{\lambda}{2c_2}\Vert\nabla\theta\Vert_{H^s(\Omega)}^2\nonumber\\
&\leq\delta\Vert\Delta q\Vert_{H^s(\Omega)}^2+C\Vert\nabla q\Vert_{H^s(\Omega)}^2+C\Vert\nabla q\Vert_{H^s(\Omega)}^6+C\Vert E\Vert_{H^s(\Omega)}^2+C\Vert E\Vert_{H^s(\Omega)}^6+C\varepsilon^2,
\end{align}
where $\delta$ is a positive constant to be determined.

 $\varepsilon\nabla q^\alpha$ and $(\ref{4.6})$ are taken as inner products, and by integrating by parts, and combining with lemma \ref{lemma2.2} and lemma \ref{lemma2.5}, we get
\begin{align}\label{4.11}
&\frac{d}{dt}\sum\limits_{|\alpha|\leq s}\int_\Omega\varepsilon\mathbf{u}^\alpha\nabla q^\alpha dx+\frac{3}{4}\Vert\Delta q\Vert_{H^s(\Omega)}^2\nonumber\\
&\leq C_1(\Vert\operatorname{div}\mathbf{u}\Vert_{H^s(\Omega)}^2+\vert\nabla\mathbf{u}\Vert_{H^s(\Omega)}^2+\Vert\nabla\theta\Vert_{H^s(\Omega)}^2)\nonumber\\
&~~~~+
C(\Vert\nabla q\Vert_{H^s(\Omega)}^2+\Vert\nabla q\Vert_{H^s(\Omega)}^6+\Vert E\Vert_{H^s(\Omega)}^2+\Vert E\Vert_{H^s(\Omega)}^6+\varepsilon^2).
\end{align}

Let $k_1=\frac{1}{4C_1}\mathrm{min}\{\mu,\mu+\nu,\frac{\lambda}{4}\}$, $k=\frac{1}{2}\mathrm{min}\{k_1,\sqrt{C_1}\}$ and $\delta=\frac{k}{4}$, we obtain the sum of $k\times(\ref{4.11})$ and $(\ref{4.10})$
\begin{align}\label{4.12}
\frac{d}{2dt}&\sum\limits_{|\alpha|\leq s}\int_\Omega(\frac{1+\varepsilon\theta^\varepsilon}{1+\varepsilon q^\varepsilon}|q^\alpha|^2+(1+\varepsilon q^\varepsilon)|\boldsymbol{u}^\alpha|^2+\frac{1+\varepsilon q^\varepsilon}{1+\varepsilon\theta^\varepsilon}|\theta^\alpha|^2+|\nabla q^\alpha|^2+2k\varepsilon\boldsymbol{u}^\alpha\nabla q^\alpha)dx\nonumber\\
&+\frac{\varepsilon\mu}{4}\Vert\nabla\boldsymbol{u}\Vert_{H^s(\Omega)}^2+\frac{\varepsilon(\mu+\nu)}{4}\Vert\mathrm{div}\boldsymbol{u}
\Vert_{H^s(\Omega)}^2+\frac{\lambda}{4c_2}\Vert\nabla\theta\Vert_{H^s(\Omega)}^2+\frac{k}{2}\Vert\Delta q\Vert_{H^s(\Omega)}^2\nonumber\\
\leq& C\Vert E\Vert_{H^s(\Omega)}^2+C\Vert E\Vert_{H^s(\Omega)}^6+C\Vert\nabla q\Vert_{H^s(\Omega)}^2+C\Vert\nabla q\Vert_{H^s(\Omega)}^6+C\varepsilon^2.
\end{align}

Define $E_1:=(E,\nabla q)$, by integrating in $[0,t], t\leq T\leq\mathrm{min}\{T^\ast,T_\varepsilon\}$ and applying Gronwall's lemma, we obtain
\begin{align}\label{4.13}
\Vert E_1(t)\Vert_{H^s(\Omega)}^2\leq C\Vert E_1(0)\Vert_{H^s(\Omega)}^2+C\int_0^t(1+\Vert E_1\Vert_{H^s(\Omega)}^4)\Vert E_1\Vert_{H^s(\Omega)}^2dt+C\varepsilon^2.
\end{align}
Let $\Lambda(t)=C\varepsilon^2\mathrm{exp}[C\int_0^t(1+\Vert E_1(\tau)\Vert_{H^s(\Omega)}^4)d\tau]$, using the Gronwall inequality, we prove that $\Lambda(t)$ is uniformly bound in $[0,T]$, for $T\in[0,\mathrm{min}\{T^\ast,T_\varepsilon\}]$.

This completes the proof of Theorem \ref{th1.3}.

\printcredits
\bibliography{cas-refs}

@article{alazard2006low,
  title={Low Mach number limit of the full {Navier-Stokes} equations},
  author={Alazard, Thomas},
  journal={Archive for rational mechanics and analysis},
  volume={180},
  pages={1--73},
  year={2006},
  publisher={Springer}
}

@article{fanelli2016highly,
  title={Highly rotating viscous compressible fluids in presence of capillarity effects},
  author={Fanelli, Francesco},
  journal={Mathematische Annalen},
  volume={366},
  pages={981--1033},
  year={2016},
  publisher={Springer}
}

@article{feireisl2012singular,
  title={A singular limit for compressible rotating fluids},
  author={Feireisl, Eduard and Gallagher, Isabelle and Novotn{\'y}, Anton{\'\i}n},
  journal={SIAM Journal on Mathematical Analysis},
  volume={44},
  number={1},
  pages={192--205},
  year={2012},
  publisher={SIAM}
}

@article{feireisl2014scale,
  title={Scale interactions in compressible rotating fluids},
  author={Feireisl, Eduard and Novotn{\'y}, Anton{\'\i}n},
  journal={Annali di Matematica Pura ed Applicata (1923-)},
  volume={193},
  number={6},
  pages={1703--1725},
  year={2014},
  publisher={Springer}
}

@article{ju2022low,
  title={Low Mach number limit of {Navier-Stokes} equations with large temperature variations in bounded domains},
  author={Ju, Qiangchang and Ou, Yaobin},
  journal={Journal de Mathematiques Pures et Appliquees},
  volume={164},
  pages={131--157},
  year={2022},
  publisher={Elsevier}
}

@article{kato1984nonlinear,
  title={Nonlinear evolution equations and the {Euler} flow},
  author={Kato, Tosio and Lai, Chi Yuen},
  journal={Journal of functional analysis},
  volume={56},
  number={1},
  pages={15--28},
  year={1984},
  publisher={Elsevier}
}

@article{korteweg1901forme,
  title={Sur la forme que prennent les {\'e}quations du mouvements des fluides si l'on tient compte des forces capillaires caus{\'e}es par des variations de densit{\'e} consid{\'e}rables mais connues et sur la th{\'e}orie de la capillarit{\'e} dans l'hypoth{\`e}se d'une variation continue de la densit{\'e}},
  author={Korteweg, Diederick Johannes},
  journal={Archives N{\'e}erlandaises des Sciences exactes et naturelles},
  volume={6},
  pages={1--24},
  year={1901}
}

@article{kwon2018incompressible,
  title={Incompressible limit of the degenerate quantum compressible {Navier-Stokes} equations with general initial data},
  author={Kwon, Young-Sam and Li, Fucai},
  journal={Journal of Differential Equations},
  volume={264},
  number={5},
  pages={3253--3284},
  year={2018},
  publisher={Elsevier}
}

@article{lions1998incompressible,
  title={Incompressible limit for a viscous compressible fluid},
  author={Lions, P-L and Masmoudi, Nader},
  journal={Journal de math{\'e}matiques pures et appliqu{\'e}es},
  volume={77},
  number={6},
  pages={585--627},
  year={1998},
  publisher={Elsevier}
}

@article{ou2023low,
  title={Low Mach number limit for non-isentropic magnetohydrodynamic equations with ill-prepared data and zero magnetic diffusivity in bounded domains},
  author={Ou, Yaobin and Yang, Lu},
  journal={arXiv preprint arXiv:2308.07745},
  year={2023}
}

@article{sha2019low,
  title={Low Mach number limit of the three-dimensional full compressible {Navier-Stokes-Korteweg} equations},
  author={Sha, Kaijian and Li, Yeping},
  journal={Zeitschrift f{\"u}r angewandte Mathematik und Physik},
  volume={70},
  pages={1-16},
  year={2019},
  publisher={Springer}
}

@article{van1894thermodynamische,
  title={Thermodynamische theorie der kapillarit{\"a}t unter voraussetzung stetiger dichte{\"a}nderung},
  author={Van der Waals, Johannes Diderik},
  journal={Zeitschrift f{\"u}r Physikalische Chemie},
  volume={13},
  number={1},
  pages={657--725},
  year={1894},
  publisher={De Gruyter (O)}
}

@article{feireisl2014inviscid,
  title={Inviscid incompressible limits on expanding domains},
  author={Feireisl, Eduard and Ne{\v{c}}asov{\'a}, {\v{S}}{\'a}rka and Sun, Yongzhong},
  journal={Nonlinearity},
  volume={27},
  number={10},
  pages={2465},
  year={2014},
  publisher={IOP Publishing}
}

@article{hattori1996existence,
  title={The existence of global solutions to a fluid dynamic model for materials for Korteweg type},
  author={Hattori, Harumi and Li, Dening},
  journal={Journal of Partial Differential Equations},
  volume={9},
  pages={323--342},
  year={1996},
  publisher={INTERNATIONAL ACADEMIC PUBLISHERS}
}

@incollection{masmoudi2007examples,
  title={Examples of singular limits in hydrodynamics},
  author={Masmoudi, Nader},
  booktitle={Handbook of differential equations: evolutionary equations},
  editor={Dafermos, Constantine and Feireisl, Eduard},
  publisher={Elsevier},
  address={Amsterdam},
  year={2007},
  pages={195--275}
}

@article{babin2001navier,
  title={{3D} {Navier-Stokes} and {Euler} equations with initial data characterized by uniformly large vorticity},
  author={Babin, A. and Mahalov, A. and Nicolaenko, B.},
  journal={Indiana University Mathematics Journal},
  year={2001},
  volume={50},
  number={1},
  pages={1--35},
  publisher={JSTOR}
}

@article{babin1999global,
  title={Global regularity of {3D} rotating {Navier-Stokes} equations for resonant domains},
  author={Babin, A. and Mahalov, A. and Nicolaenko, B.},
  journal={Indiana University Mathematics Journal},
  year={1999},
  volume={48},
  number={3},
  pages={1133--1176},
  publisher={JSTOR}
}

@article{clopeau1998vanishing,
  title = {On the vanishing viscosity limit for the {2D} incompressible {Navier-Stokes} equations with the friction type boundary conditions},
  author = {Clopeau, Thierry and Mikelic, Andro and Robert, Raoul},
  journal = {Nonlinearity},
  year = {1998},
  volume = {11},
  number = {6},
  pages = {1625--1636},
  publisher = {IOP Publishing}
}

@article{BourguignonBrezis1974,
    author = {Bourguignon, J.-P. and Brezis, H.},
    title = {Remarks on the {Euler} equation},
    journal = {Journal of Functional Analysis},
    volume = {15},
    year = {1974},
    pages = {341--363},
    issn = {0022-1236},
    doi = {10.1016/0022-1236(74)90027-5},
}

@article{hou2018global,
  title={Global well-posedness of the {3D} non-isothermal compressible fluid model of Korteweg type},
  author={Hou, Xiaofeng and Peng, Hongyun and Zhu, Changjiang},
  journal={Nonlinear Analysis: Real World Applications},
  volume={43},
  pages={18--53},
  year={2018},
  publisher={Elsevier}
}

@article{hou2017vanishing,
  title={Vanishing capillarity limit of the compressible non-isentropic {Navier-Stokes-Korteweg} system to {Navier-Stokes} system},
  author={Hou, Xiaofeng and Yao, Lei and Zhu, Changjiang},
  journal={Journal of Mathematical Analysis and Applications},
  volume={448},
  number={1},
  pages={421--446},
  year={2017},
  publisher={Elsevier}
}

@article{ou2022incompressible,
  title={Incompressible Limit of Isentropic Navier--Stokes Equations with Ill-Prepared Data in Bounded Domains},
  author={Ou, Yaobin and Yang, Lu},
  journal={SIAM Journal on Mathematical Analysis},
  volume={54},
  number={3},
  pages={2948--2989},
  year={2022},
  publisher={SIAM}
}
\end{document}